%%%%%%%%%%%%%
% Andrei Tyurin (Steklov Inst., Moscow and Univ. of Warwick)
% Special Langrangian geometry and slightly deformed algebraic geometry
% (spLag and sdAG)
% Uses: amstex, amsppt 2.1 (and epsf.tex for 2 inessential figures)
%%%%%%%%%%%%%

 \input amstex
 \documentstyle{amsppt}

 \input epsf % comment out if you don't want the figures
 \epsfverbosetrue

 \magnification=\magstep1
 \vsize=24.2true cm
 \hsize=15.3true cm
 \nopagenumbers\topskip=1truecm
 \headline={\tenrm\hfil\folio\hfil}

 \TagsOnRight

 \define\Span#1{\left< #1 \right>}
 \define\half{\textstyle\frac12}
 \define\1{^{-1}}
 \define\wave{\widetilde}
 \define\tensor{\otimes}
 \define\dd{\roman{d}}
 \define\id{\operatorname{id}}
 \define\pr{\operatorname{pr}}
 \define\Tr{\operatorname{Tr}}
 \redefine\det{\operatorname{det}}
 \define\II{\roman{II}}
 \define\can{\roman{can}}
 \define\suniv{_{\roman{univ}}}
 \define\irr{\roman{irr}}
 \define\red{\roman{red}}
 \define\LC{\roman{LC}}
% \define\red{_{\roman{red}}} % reduced subscheme
% \define\mi{_{\roman{min}}} % minimal surface
% \define\rest#1{_{\textstyle{|}#1}} % restriction of map to subset

% Script letters
% \define\sA{{\Cal A}} % sheaf of algebras A
% \define\sB{{\Cal B}} % conductor ideal
% \define\sG{{\Cal G}} % sheaf G
% \define\sJ{{\Cal J}} % sheaf J
 \define\sI{{\Cal I}}
 \define\sL{{\Cal L}}
 \define\sM{{\Cal M}}
 \define\Oh{{\Cal O}}
% \define\sS{{\Cal S}}
% \define\sF{{\Cal F}}
% \define\sE{{\Cal E}}

% short Greeks
 \define\al{\alpha}
 \define\be{\beta}
 \define\de{\delta}
 \define\ep{\varepsilon}
 \define\fie{\varphi}
% \define\ga{\gamma}
 \define\ka{\kappa}
 \define\om{\omega}
 \define\si{\sigma}
% \define\De{\Delta}
 \define\Ga{\Gamma}
 \define\La{\Lambda}
 \define\Om{\Omega}
% \define\la{\lambda}
 \define\Si{\Sigma} 

 \define\proj{{\Bbb P}}
 \define\C{{\Bbb C}}
 \define\Q{{\Bbb Q}}
 \define\R{{\Bbb R}}
 \define\Z{{\Bbb Z}}

% \mathops
 \define\ad{\operatorname{ad}} % adjoint bundle
 \define\adj{\operatorname{adj}}
% \define\alg{\operatorname{alg}}
% \define\res{\operatorname{res}}
 \define\codim{\operatorname{codim}} % codimension
% \define\diff{\operatorname{diff}}
% \define\topo{\operatorname{top}}
% \define\var{\operatorname{var}}
% \define\vcodim{\operatorname{v.codim}}
% \define\vdim{\operatorname{v.dim}}
% \define\Diff{\operatorname{Diff}}
% \define\Ext{\operatorname{Ext}}
% \define\Grass{\operatorname{Grass}} % Grassmann variety
% \define\Mod{\operatorname{Mod}}
% \define\Mon{\operatorname{Mon}}
 \define\Pic{\operatorname{Pic}} % Picard scheme
 \define\Aut{\operatorname{Aut}}
 \define\rk{\operatorname{rank}}
 \define\im{\operatorname{im}}
 \define\tr{\operatorname{tr}}
 \define\End{\operatorname{End}}
 
 \define\LaGr{\operatorname{\La_{\uparrow}\!}}
 \define\orGr{\operatorname{Gr_{\uparrow}\!}}
 \define\Hom{\operatorname{Hom}}
 \define\Mod{\operatorname{Mod}}
 \define\Rep{\operatorname{Rep}}
 \redefine\Re{\operatorname{Re}}
 \define\Sing{\operatorname{Sing}}
 \define\Spec{\operatorname{Spec}}
 \define\Vol{\operatorname{Vol}}
% \define\Stab{\operatorname{Stab}} % Stabiliser group
% \define\Sym{\operatorname{Sym}}

% Lie groups
% \redefine\O{\operatorname{O}}
 \define\SO{\operatorname{SO}}
 \define\U{\operatorname U} % unitary group
 \define\PU{\operatorname{PU}}
 \define\SU{\operatorname{SU}}
 \define\sSU{\operatorname{{\Cal{SU}}}}
 \define\su{\operatorname{{\frak{su}}}}
 \define\Spin{\operatorname{Spin}}
% \define\SP{\operatorname{Sp}}
% \define\odd{\roman{odd}}
% \define\pt{\roman{pt}}

 \document\topmatter
 \title Special Langrangian geometry and slightly deformed algebraic
 geometry (spLag and sdAG)
 \endtitle
 \author Andrei Tyurin \endauthor

 \address Algebra Section, Steklov Math Institute, Ul\. Gubkina 8,
 Moscow, GSP--1,
 117966, Russia \endaddress
 \email Tyurin\@tyurin.mian.su
 {\it or}\newline Tyurin\@Maths.Warwick.Ac.UK
 {\it or}\newline Tyurin\@mpim-bonn.mpg.de
 \endemail

 \abstract
The special geometry of calibrated cycles, closely related to mirror
symmetry among Calabi--Yau 3-folds, is itself a real form of a new
subject, which we call {\it slightly deformed algebraic geometry}. On the
other hand, both of these geometries are parallel to classical gauge
theories and their complexifications. This article explains this
parallelism, so that the appearance of invariants of new type in
complexified gauge theory (see [D-T] and [T]) can be accompanied by
analogous invariants in the theory of special Lagrangian cycles, for which
the development is at present much more modest than in gauge theory. We
discuss related geometric constructions, arising from mirror symmetry and
symplectic geometry.

 \endabstract

 \endtopmatter
 \rightheadtext{spLag and sdAG geometries}

 \head \S1. spLag cycles \endhead

We begin by recalling the geometric construction for a pair $\sL\subset
S$, where $S$ is a smooth symplectic manifold of dimension $2n$ with a
given tame almost complex structure $I$, and $\sL\subset S$ a smooth,
oriented Lagrangian submanifold (of maximal dimension
$\dim\sL=n=\half\dim S$); this is now quite popular in the set-up of
Calabi--Yau threefolds. The structure on $S$ is an {\it almost K\"ahler
structure}, and we say for short that $S$ is an {\it aK manifold}. Write
$\om$ for the symplectic form and $I$ for the almost complex structure on
$S$, so that the tangent space $TS_p$ at a point $p$ is $\C^n$ with the
constant symplectic form $\Span{\ ,\,}=\om_p$ and the constant Euclidean
metric $g_p$, giving the Hermitian triple $(\om_p,I_p,g_p)$.

We now define the {\it Lagrangian Grassmannian} $\LaGr p=\LaGr(TS_p)$ to
be the Grassmannian of maximal oriented Lagrangian subspaces in $TS_p$.
Taking this space over every point of $S$ gives the {\it oriented
Lagrangian Grassmannisation} of $TS$
 $$
 \pi \colon \LaGr(S)\to S \quad\text{with} \quad \pi\1(p)=\LaGr p.
 \tag1.1
 $$

Our tame almost complex structure on $S$ gives each fibre the standard
form
 $$
 \LaGr p=\U(n) / \SO(n)
 \tag1.2
 $$
This space admits a canonical map
 $$
 \det\colon\LaGr p\to\U(1)=S^1_p
 \quad\text{sending $u\in\U(n)$ to $\det u\in\U(1)=S^1$.}
 \tag1.3
 $$
Recall that the inverse image of the fundamental class of $S^1$ on $\LaGr
p$ is the {\it universal Maslov class}. Taking this map over every point
of $S$ gives the map 
 $$ 
 \det \colon \LaGr S\to S^1(L_{-K}),
 \tag1.4
 $$
where $S^1(L_{-K})$ is the unit circle bundle of the line bundle
$\bigwedge^n TS=\det TS$, with first Chern class
 $$
 c_1(\det TS)=-K_S,
 \tag1.5
 $$
where $K_S$ is the canonical class of $S$. Recall that, as a cohomology
class, $K_S$ does not depend on the compatible almost complex structure.

Now for every oriented Lagrangian cycle $\sL\subset S$, we have
the Gauss lift of the embedding $i\colon \sL\to S$ to a section 
 $$
 G(i) \colon \sL\to \LaGr(S) \vert_{\sL}
 \tag1.6
 $$ sending a point $p\in \sL$
to the oriented subspace $T\sL_p\subset TS_p$. The composite of
this Gauss map with the projection (1.4) gives the map 
 $$
 {\det} \circ G(i)\colon \sL\to S^1(L_{-K}) \vert_{\sL} 
 \tag1.7
 $$

Now suppose that the cohomology class of the symplectic form is
proportional to the canonical class of $S$, that is,
 $$ 
 \ka \cdot[\om]=K_S \quad\text{for some $\ka\in \Q$};
 \tag1.8
 $$
then the restriction $\det TS\vert_{\sL}$ is topologically trivial,
because the restriction of $[\om]$ to a Lagrangian $\sL$ is zero.
Moreover, suppose that $K_S$ has a Hermitian connection $a_K$ with
curvature form
 $$
 F_{a_K}=(2\pi i) \om,
 $$ 
where $\om$ is our symplectic form. Then this connection restricts to a
flat connection on $\sL$. Let
 $$
 \Span{\sL}
 =\bigl\{B\subset S\bigm|\sL\subset B,\pi_1(B)=1\text{ and } F_a
 \vert_B=0\bigr\}
 $$
be the maximal simply connected submanifold containing $\sL$ on which
$K_S$ and $a$ restrict to a trivial bundle and a flat connection.

Then on $\Span{\sL}$ there exists a canonical trivialisation 
 $$ 
 S^1(L_{-K}) \vert_{\Span{\sL}}=\Span{\sL} \times S^1 
 \tag1.9 
 $$
that preserves the Hermitian form and the canonical projection 
 $$ 
 \pr \colon S^1(L_{-K}) \vert_{\sL}\to S^1.
 \tag1.10 
 $$
Now composing the maps (1.4), (1.6) and (1.10) gives a map 
 $$
 m={\pr}\circ{\det}\circ G(i)\colon\sL\to S^1
 \tag1.11 
 $$

\definition{Definition 1.1} (1) $m$ is called the {\it phase map};

(2) $\sL$ is called a {\it special Lagrangian cycle} of $S$ ({\it spLag
cycle} for short) if $m(\sL)$ is a point or, equivalently, the
differential of $m$ vanishes:
 $$
 \dd m=0.
 \tag1.12
 $$
\enddefinition

 \remark{Remark} In this definition, we call $\sL$ a cycle rather than a
submanifold, because it may be singular. We really only need the Gauss
map (1.6) to be well defined; thus $\sL$ can have nodes, and so on. Thus
below we call a {\it cycle} any subvariety with regular Gauss map.
 \endremark

 \remark{Mirror digression} (1.8) holds automatically if $S$ is a
Calabi--Yau manifold, that is, $K_S=0$. In other words, for any
Lagrangian submanifold $\sL$
 $$
 \Span{\sL}=S,
 \tag 1.13
 $$
and we have the map $m\colon\sL\to S^1$ (1.11). In this case, the notion
of spLag cycle coincides with that in calibrated geometry (see [H-L]).

Recall that a {\it complex orientation} of a Calabi--Yau manifold $X$ is a
choice of trivialisation of the canonical line bundle $L_K$, that is, a
holomorphic 3-form $\theta$. For an oriented Calabi--Yau threefold
$(X,\theta)$, a spLag cycle is a 3-dimensional Lagrangian submanifold $\sL$
such that the restriction 
 $$
 \Re\theta \vert_{\sL}=0.
 \tag 1.14
 $$
 \endremark
The local deformation theory of such submanifolds is well understood. The
tangent space to the moduli space $\sM_{\sL}$ of such deformations at a
submanifold $\sL$ is $H^1(\sL,\R)$, viewed as the space of harmonic 1-forms
on $\sL$. This space doesn't depend on the second quadratic form or others
attributes of embeddings. In particular, if $H^1(\sL,\R)=0$ then $\sL$
is rigid as a special Lagrangian submanifold. So we can expect that there
exists a finite set of such submanifolds $\{\sL_1,\dots, \sL_N\}$ in
one cohomology class $[\sL_i]$. This subject is quite popular now, and we
would like to remark that this construction also works for Fano varieties.

 \remark{Chern--Simons digression} Let $\Si$ be a compact, smooth,
oriented Riemann surface of genus $g$, and
 $$
 \pi_1(g)=\Span{a_1, \dots, a_g, b_1, \dots, b_g
 \bigm| \Pi_{i=1}^{g} [a_i, b_i]=\id}
 \tag1.15
 $$
the usual presentation of the fundamental group $\pi_1(g)$.

The space $\Rep\pi_1(g)$ is the moduli space of classes of $\SU(2)$
representations of $\pi_1(g)$, and is smooth as an orbifold. But, as a
manifold, it is singular at the reducible representations
 $$
 \Sing\Rep\pi_1(g)=\Rep\pi_1(g)^{\red}.
 \tag1.16
 $$
The space $\Rep\pi_1(g)$ contains the special subspace of representations
that are trivial on the $b_i$
 $$
 B=\bigl\{\rho\in\Rep\pi_1(g)\bigm|\rho(b_i)=\id
 \text{ for } i=1,\dots,g\bigr\}.
 \tag1.17
 $$
As usual,
 $$
B^{\irr}\subset\Rep\pi_1(g)^{\irr}
 \tag1.18
 $$
is the subset of irreducible representations.

Then $B$ is a Lagrangian suborbifold of $\Rep\pi_1(g)$ with respect to
the canonical symplectic form $\om_{\Si}$ (1.26). To apply the phase
map to this Lagrangian submanifold $\sL=B$ in $S=R_g$ we have to
remark that
 $$
[\om_{\Si}]=-4 \cdot K_{R_g}
 \tag1.19
 $$
(see [N-R], [R]).

Now fixing a complex structure on $\Si$ gives a complex structure on
$R_g$, the Weil--Peterson metric on $R_g$, which is K\"ahlerian with
symplectic form $\om_{\Si}$ and a unitary connection on the anticanonical
bundle with curvature form $i\cdot\pi\om_\Si$ and the phase map 
 $$
 m_g\colon B^{\irr}\to S^1=\U(1) 
 \tag1.20
 $$
which can be investigated in the usual way (as in [A-M]).

Bogomolov remarked that there exists a complex structure on $\Si$ for
which $B$ is a spLag cycle. Indeed, consider the complex structure on
$\Si$ of genus 2 such that all the cycles $b_i$ are {\it real}. For example if
$\Si$ is a double of a pair of pants:
 $$
 \epsfbox{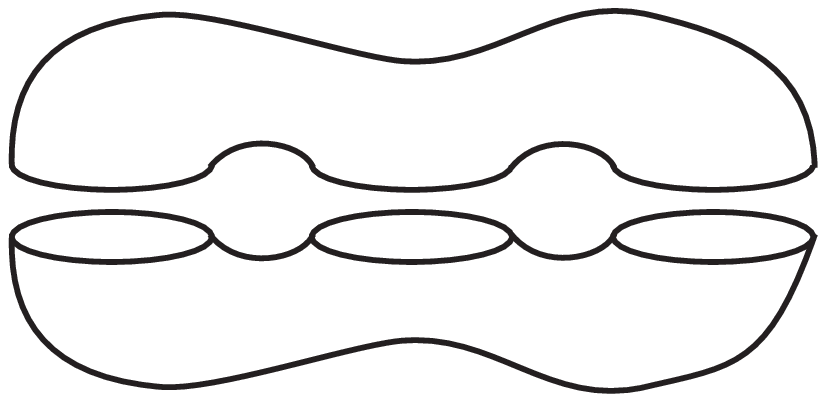}
 $$
we get a compact complex Riemann surface of genus 2 with a real structure
such that cycles $b_1$ and $b_2$ are real Then the usual argument implies
that the cycle $B$ must be real. Then one can see that $B$ is spLag (see
1.14).

It is easy to see that $B$ is rigid as a spLag cycle. The following
question is very interesting:

 \remark{Question} Is there another spLag cycle $B'$ in the same homology
class $[B]=[B']$?
 \endremark

Now there is an absolutely canonical almost K\"ahler structure on $R_g$,
constructed by Guruprasad and Nilakantan. To describe it, we use the
standard presentation (1.15) of $\pi_1(g)$ and another ``dual''
presentation given by the following construction (see [G-N]): let 
 $$
 r_i=\Pi_{j=1}^{i} [a_j, b_j], \quad \text{and set} \quad
 \al_i=r_{i-1} b_i\1r_i\1\quad \text{and} \quad
 \be_i=r_i a_i\1r_{i-1}\1.
 \tag1.21
 $$
Then
 $$
 \pi_1(g)=\Span{\al_1, \dots, \al_g, \be_1, \dots,\be_g \bigm|
 \Pi_{j=1}^{g} [a_j, b_j]=1}
 \tag$1.15'$
 $$
is another presentation of $\pi_1(g)$. Sending the generators $a_i, b_j$
to $\al_i,\be_j$ gives an involutive automorphism $W$ of $\pi_1(g)$, that
is, and element $W\in\Mod_g$ with $W^2=\id$. 
 
To describe the tangent space $(TR_g)_{\rho}$, recall that a $\SU(2)$
representation $\rho$ makes the Lie algebra $\su(2)$ into a
$\pi_1(g)$-module $\su(2)_{\rho}$ by the adjoint action, and,
as the tangent space to an orbifold point, 
 $$
 (TR_g)_{[\rho]}=H^1(\pi_1(g),\su(2)_{\rho}).
 \tag1.22
 $$
The group of cycles of the module $\su(2)_{\rho}$ is the group of skew
homomorphisms from $\pi_1(g)$ to this module, that is
 $$
 Z^1(\su(2)_{\rho})=
 \bigl\{u \colon\pi_1(g)\to \su(2)_{\rho} \bigm|
 u(g_1\circ g_2)=u(g_1)+g_1 (u(g_2)) \bigr\}.
 \tag1.23
 $$
Of course, any such $u$ extends to a $\Z$-linear map of the integral
group algebra:
 $$
 u \colon\Z\pi_1(g)\to \su(2),
 \tag1.24
 $$
and the boundary subspace is
 $$
 B^1(\su(2)_{\rho})=
 \bigl\{u \colon\pi_1(g)\to \su(2)_{\rho} \bigm| u(g)=g(v)-v
 \text{ for some } v\in \su(2)\}.
 \tag1.25
 $$
Then
 $$
 (TR_g)_{[\rho]}=Z^1(\su(2)_{\rho}) / B^1(\su(2)_{\rho}).
 $$

The canonical symplectic structure $\om_\Si$ on $R_g$ can be defined as
a skewsymmetric bilinear form on the tangent space of every class of
representations $[\rho]$ by the $\PU(2)$ invariant bilinear form on
$Z^1(\su(2)_{\rho})$ given by 
 $$
 \Span{u,v}=\sum_{i=1}^g
 \Bigl[(u(r_{i-1}\1-b_i\1\cdot r_i\1),
 v(a_i))+(u(a_i\1r_{i-1}\1-r_i\1), v (b_i)) \Bigr]
 \tag1.26 
 $$ 
where $(\ ,\,)$ is the standard inner product $(m,m)=-\tr m^2$ on
$\su(2)$.

The inner product on the space of cycles $Z^1(\su(2)_{\rho})$ (1.23) is
given by the formula
 $$
 G (u, v)=\sum_{i=1}^g \Bigl[(u(\al_i), v (\al_i))+(u(\be_i),
 v (\be_i))\Bigr]
 \tag1.27
 $$
(see [G-N])).

 \proclaim{Proposition 1.1 \rm [G-N]} 1) The inner product (1.27) is
positive definite.

2) The restriction of the inner product (1.27) to the orthogonal
 $$
 B^1(\su(2)_{\rho})^{\perp}=(TR_g)_{[\rho]}
 \tag1.28
 $$
defines a special Riemannian metric on $R_g^{\irr}$ and on $R_g$ as an
orbifold. \endproclaim
 
 \definition{Definition 1.2} This metric is called the {\it Guruprasad
metric} ({\it GN metric} for short). \enddefinition

It can be checked (see [G-N], Proposition 3.1) that the GN metric is
compatible with the symplectic form (1.26); this canonical Hermitian
structure and the canonical connection on the tangent bundle (the
Levi--Civita connection) induce a Hermitian structure and a connection on
the determinant line bundle $L_{-K}$ whose curvature coincides with the
symplectic form. Thus the phase map
 $$
 m_{G}\colon B\to S^1
 \tag1.29
 $$
is absolutely canonical.
 \endremark

Now we return to the general phase map.

There is a standard way of describing the differential of the phase map:
let
 $$
\nabla_\LC \colon \Ga (TS)\to \Ga(TS\tensor T^*S)
 \tag1.30
 $$
be the Levi--Civita connection, which we restrict to the restriction of
the tangent bundle to a Lagrangian cycle $\sL$. Let $N_{\sL\subset S}$
be the normal bundle of $\sL$ in $S$. Then the Levi--Civita connection
(1.30) defines connections 
 $$
 \aligned
\nabla_\LC \colon \Ga (T \sL) &\to \Ga(T \sL\tensor T^* \sL)\\
\nabla_\LC \colon \Ga (N_{\sL\subset S}) &\to \Ga(N_{\sL\subset S}\tensor T^* \sL)
 \endaligned
 $$
and the tensor
 $$
\II\colon T \sL\to \Hom(T \sL, N_{\sL\subset S})
 $$
(called the second quadratic form), that is,
 $$
\II\in\End T \sL\tensor N_{\sL\subset S}.
 \tag1.31
 $$
The trace component of $\End T \sL=\C \oplus \ad T \sL$ gives
 a section
 $$
H\in \Ga (N_{\sL\subset S}),
 \tag1.32
 $$
called the mean curvature section.

Let $V$ be a vector field on $\sL$. Then pointwise, the value of the
differential $\dd m_G$ on $V$ is given by the inner product
 $$
 \dd m_G(V_p)=(H_p, I_G(V_p)),
 \tag1.33
 $$
where $I_G$ is the operator of the Guruprasad almost complex structure.

\head \S2 sdAG cycles (slightly deformed Algebraic Geometric cycles)
\endhead

The construction of the phase map can be ``complexified'' to give the
complex phase map 
 $$
 m_\C \colon \Si\to \C^*
 \tag2.1
 $$
for any middle dimensional cycle $\Si\subset S$, where $S$ is an almost
K\"ahler manifold.

We illustrate this construction in the case of an aK surface $S$ (see the
beginning of \S1). Then over any point $p\in S$ we have as the fibre of
the tangent bundle $TS_p=\C^2$ with the standard Hermitian triple
$(\om_p,I_p,g_p)$, and the oriented Lagrangian Grassmannian is a subspace
of the oriented Grassmannian
 $$
 \LaGr p\subset\orGr(2, TS_p).
 \tag2.2
 $$
Now we saw that
 $$
\LaGr p=S^2 \times S^1
 \tag2.3
 $$
(see (1.2)), and
 $$\orGr(2, TS_p)=S^2_- \times S^2_+,
 \tag2.4
 $$
so that what we need to see is that the inclusion (2.2) is componentwise,
that is,
 $$
S^2=S^2_- \quad \text{and} \quad S^1\subset S^2_+
 \tag2.5
 $$
and the first identification and the second inclusion carry deep geometric
meaning in the theory of Riemannian twistors.

First of all, every oriented plane $V\subset TS_p$ defines a new complex
structure $I_V$ on $TS_p$ which is compatible with the Euclidean metric
$g_p$. Thus we have the natural map
 $$
 I \colon\orGr(2, TS_p)\to S^2=\proj W^+_p,
 \tag2.6
 $$
where $W^+=\C^2$ is the positive spinor space of the Euclidean metric.

The fibre of this map 
 $$
I\1(I_V)=\proj \C^2_V=\proj W^-_p
 \tag2.7
 $$
is the $I_V$ projective line -- that is, complex directions in the complex
plane with complex structure $I_V$. It is easy to see that the
decomposition (2.4) is precisely the twistor decomposition 
 $$
 \orGr(2, TS_p)=\proj W^-_p \times \proj W^+_p.
 \tag2.8
 $$
Now for any oriented plane $V$, consider its {\it K\l
ahler angle}
$\al(V)$, given by the formula
 $$
 \om_p \vert_V=\arccos\al \cdot\Vol_{g_p}.
 \tag 2.9
 $$
The geometric meaning of $\al (V)$ is following: it is the angle between
$V$ and $I_p(V)$. Sending an oriented plane to the K\"ahler angle gives a
map
 $$
 \al_p \colon\orGr(2, TS_p)\to [0,\pi]
 \tag2.10
 $$
of our Grassmannian to the interval.

The angle $\al (V)$ only depends on the complex structure $I_V$, so that
the map $\al_p$ is the composite
 $$
 \orGr(2, TS_p) @> \pr_+>> \proj W^+ @>h>> [0,\pi]
 \tag2.11
 $$
and the last map
 $$
h \colon \proj W^+=S^2\to [0,\pi]
 \tag2.12
 $$
is the standard {\it height function} on $S^2$ sending the standard
sphere in Euclidean $\R^3$ to the third coordinate:
 $$
 \epsfbox{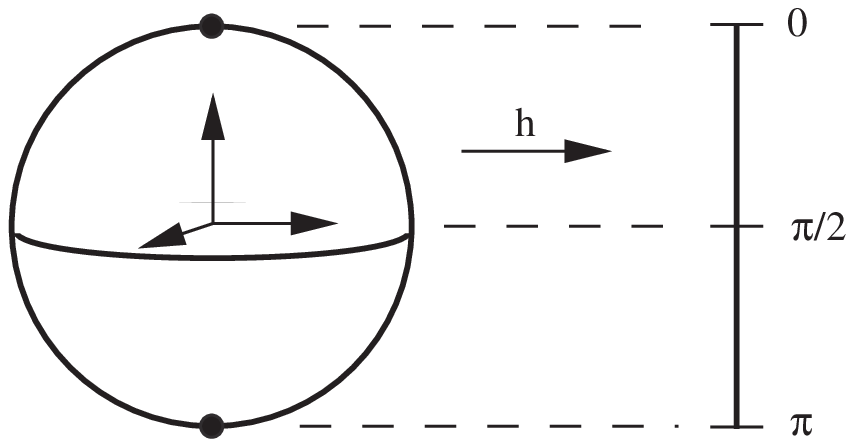}
 $$
The fibres of the K\"ahler angle map $\al$ have the following geometric
interpretration:
 $$
 \al\1(0)=\proj\C^2=\proj W^-
 \tag2.13
 $$
is the space of complex tangent directions at a point $p$. In the same
vein,
 $$
 \al\1(\pi)=\overline{\proj W^-}
 \tag2.14
 $$
is the space of complex directions at a point $p$, but with the opposite
orientation. Now
 $$
 \al\1(\pi/2)=\LaGr p
 \tag2.15
 $$
is the Lagrangian Grassmannian (1.1). It is easy to see that the
involution changing orientation sends $\al$ to $\pi-\al$ and preserves the
oriented Lagrangian Grassmannian (2.15). So the decomposition of the
Lagrangian Grassmannian (2.3) admits the following geometric meaning:
 $$
 \LaGr p=\proj W^- \times h\1(\pi/2).
 \tag2.16
 $$
So we can see that, in the same way that the fibres of the height
map $h$ are the same ($=S^1$) except for two points, the north and south
poles, the fibres of the K\"ahler angle map are the same except for the
two spheres
 $$
 (\proj W^- \times 0) \cup (\proj W^- \times\pi)
 \tag2.17
 $$
Thus, removing these exceptional fibres, we get the map
 $$
 \det_\C \colon\orGr (2, TS_p) \setminus
 (\al\1(0)\cup\al\1(\pi)) @>\pr_+>> \C^*,
 \tag 2.18
 $$
whose restriction to the Lagrangian Grassmannian $\LaGr p$ gives the
determinant map $\det$ of (1.3):
 $$
 \det_\C \vert_{\LaGr p}=\det.
 \tag2.19
 $$
Moreover, the image of the last map is the circle
 $$
 S^1=\al_p\1(\pi/2)\subset \C^*\subset \proj W^+.
 \tag2.20
 $$

Now we can give another construction for the map $\det_\C$: we fix any
nonzero element $w\in\bigwedge^2\C^2$, and call it a {\it complex
orientation} of $\C^2$. Let us fix a standard orthogonal basis
$e_1,e_2,I(e_1),I(e_2)$ in $\C^2$ such that
 $$
w=(e_1+I(e_1))\wedge (e_2+I(e_2)).
 \tag2.21
 $$
Now for any oriented plane $V\in \orGr (2, TS_p)$ consider the oriented
orthogonal basis $v_1, v_2$ and let $M_V$ be the linear transformation
sending $e_i$ to $v_i$ and $I(e_i)$ to $I(v_i)$. It is easy to check the
following result: (see for example [H-L])

 \proclaim{Lemma 2.1} 1) $M_V$ is $\C$-linear trasformation $TS_p\to
TS_p$ that is
 $$
 M_V\in\End_\C (\C^2)
 \tag2.22
 $$
is a complex matrix;

2) $\det_\C=0$ iff $V$ is a complex or an anticomplex direction
and 
 $$
V=\ker M_V;
 \tag2.23
 $$

3) $V$ is Lagrangian iff $M_V\in\U(2)$ is unitary matrix.

 \endproclaim

So sending $V$ to $\det M_V$ we get the map
 $$
d\colon \orGr (2, TS_p)\to \C,
 \tag2.24
 $$
and it is easy to see the following:

 \proclaim{Lemma 2.2} 1) 
 $$
 \vert d (V) \vert \leq 1
 \tag2.25
 $$
and
 $$
 \vert d (V) \vert=1 \iff V\in \LaGr p;
 \tag2.26
 $$

2) let 
 $$
V\to \overline{V}
 \tag2.27
 $$
be the change of orientation map then
 $$
 d (V)=- d (\overline{V}).
 \tag2.28
 $$
 \endproclaim
 Thus the map
 $$
 d^2 \colon \orGr (2, TS_p)\to \C
 \tag2.29
 $$
is the {\it orientation-forgetting} map.

Now consider the decomposition
 $$
 \orGr (2, TS_p)=\orGr (2, TS_p)^+ \cup d\1(S^1) \cup \orGr (2, TS_p)^- 
 \tag2.30
 $$
where
 $$
\orGr (2, TS_p)^-=\overline{\orGr (2, TS_p)^+}.
 \tag2.31
 $$
Then we have two maps
 $$
 \aligned
 d \colon \orGr (2, TS_p)^+ &\to \{\vert z \vert \leq 1\}\subset \C \\
 - (\overline{d})\1\colon \orGr (2, TS_p)^- &\to
 \{\vert z \vert \ge 1\}
 \endaligned
 \tag2.32
 $$
both of whose restrictions to $d\1(S^1)$ coincide:
 $$
 d \vert_{d\1(S^1)}=(\overline{d})\1\vert_{d\1(S^1)}.
 \tag2.33
 $$
Thus we get:

 \proclaim{Lemma 2.3} 1) These two maps glue to a map
 $$
 \det_\C \colon \orGr (2, TS_p)\to \proj^1=S^2.
 \tag2.34
 $$

2) This map is {\it smooth} and equal to the map $\pr_+$ of (2.11).
 \endproclaim

Recall that by the theory of Riemannian twistors, the real structure
acts on the projective twistor space $\proj W^+$ as the antipodal map
$z\to -(\overline{z})\1$.

\subheading{Globalisation} Now we have to globalise these constructions,
just as we did with the oriented Lagrangian (1.1). First of all, we have
the Grassmannisation of the tangent bundle of $S$
 $$
 \orGr(2,TS)=\proj W^- \times_S \proj W^+,
 \tag2.35
 $$
the bundle of complex quadrics on $S$ which fibrewise is the product of
twistor spaces. Of course this bundle contains the Lagrangian
Grassmannisation (1.1)
 $$
 \LaGr\subset\orGr(2,TS),
 \tag2.36
 $$
and there exists a natural projection
 $$
 \orGr(2, TS) @> \pr_+>> \proj W^+.
 \tag2.37
 $$
The K\"ahler angle map
 $$
 \al \colon\orGr (2, TS)\to [0,\pi] 
 \tag2.38
 $$
(which fibrewise is the map (2.10)) factors as the composite
 $$
 \orGr (2, TS) @> \pr_+>> \proj W^+ @>\al>> [0,\pi]
 \tag2.39
 $$
and
 $$
 \LaGr=\al\1(\pi/2).
 \tag2.40
 $$
Now the spinor bundle $\proj W^+$ is the projectivisation 
 $$
 \proj W^+=\proj(\Oh_S \oplus \Oh_S(-K_S)).
 \tag2.41
 $$
In particular, this projective bundle contains the unitary circle bundle
of the anticanonical line bundle
 $$
 S^1(L_{-K})\subset \proj W^+
 \tag2.42
 $$
and the restriction to $\LaGr$ of the $\pr_+$ map of (2.39) gives the
map (1.4).

Now for any oriented 2-dimensional submanifold $\Si\subset S$ one can
consider the Gauss lift of the embedding $i\colon \Si\to S$ to the
section
 $$
 G(i) \colon \Si\to\orGr (2, TS) \vert_\Si=\proj W^- \times_S
 \proj W^+ \vert_\Si,
 \tag2.43
 $$
sending a point $p\in\Si$ to the oriented subspace $T\Si_p\subset TS_p$.
The composite of this Gauss map with the projection $\pr_+$ (2.11)
defines the map
 $$
 \det_\C \colon \Si\to \proj W^+ \vert_\Si.
 \tag2.44
 $$
Now let $a_\LC$ be the Levi--Civita connection on the line bundle
$L_{-K}$ and 
 $$
 F_\LC\in \Om^2(S)
 \tag2.45
 $$
its curvature, viewed as a 2-form on $S$.

\definition{Definition 2.1} A cycle $\Si$ is called {\it canonically flat}
if there exists a simply connected submanifold $B \supset \Si$ such that
 $$
 F_\LC \vert_{B}=0.
 \tag2.46
 $$
\enddefinition

If $\Si$ is canonically flat, then there exists a canonical trivialisation
 $$
 \proj W^+ \vert_\Si=\Si \times S^2 
 \tag2.47
 $$
and a canonical projection
 $$
 \pr \colon \proj W^+ \vert_\Si\to S^2=\proj^1.
 \tag2.48
 $$
Moreover, the composite of all the maps gives the map
 $$
 m_\C={\pr} \circ {\det_\C} \circ G(i) \colon \Si\to S^2=\proj^1
 \tag2.49
 $$
\definition{Definition 2.2} 1) This map is called the {\it complex phase
map};

2) $\Si$ is called a sdAG cycle if $m_\C(\Si)$ is a point, or equivalently
the differential
 $$
 \dd m_\C=0.
 \tag2.50
 $$
\enddefinition

In particular, a spLag cycle is a sdAG cycle. We explain below what
``sdAG'' stands for.

 \remark{Mirror digression} As before, (2.46) holds automatically if $S$ is
a simply connected Calabi--Yau surface, that is, a K3 surface. In this
case, again $B=S$. The target space $S^2=\proj^1$ of the complex phase map
(2.49) is the space of integrable complex structures compatible with the
Calabi--Yau metric $g$ on $S$. Thus any sdAG cycle $\Si$ is a complex
curve for the complex structure $m_\C (\Si)\in S^2=\proj^1$ and any spLag
cycle is a complex curve for the specific complex structure on $S$. Any
such curve must be a fibre of an elliptic pencil or a $-2$-rational curve.
 \endremark

The target space $S^2$ admits the standard complex structure
 $$
S^2=\C \cup\infty=D^+ \cup S^1 \cup D^-,
 \tag2.51
 $$
where
 $$
 D^+=\bigl\{z \bigm| \vert z \vert<1 \bigr\};
 \quad D^-=\bigl\{z \bigm| \vert z \vert>1 \bigr\};
 \quad S^1=\bigl\{z \bigm| \vert z \vert=1 \bigr\}.
 \tag2.52
 $$
For any oriented cycle $\Si$ we have the following cases:
 $$
 \align
 m_\C (\Si)\subset D^+ &\iff \Si \quad \text{is symplectic};\\
 m_\C (\Si)\subset S^1 &\iff \Si \quad \text{is Lagrangian};\\
 m_\C (\Si)=0 &\iff \Si \quad \text{is complex or algebraic,
AG cycle for short};\\
 m_\C (\Si)=\infty &\iff \Si \quad \text{is anticomplex, antiAG
cycle for short};
 \endalign
 $$
Thus {\it sdAG} is an acronym for {\it slightly deformed Algebraic
Geometric} cycle: every point of the target sphere of the complex phase
map can be deformed along its meridian to 0 -- the image of the complex
phase map of complex cycles.

Suppose $\Si$ is a symplectic oriented cycle in $S$. Then the image
$m_\C(\Si)\subset D^+$ is compact, and there is a unique minimal disc
containing $m_\C(\Si)$
 $$
 m_\C (\Si)\subset D_\Si\subset D^+.
 \tag2.53
 $$

\definition{Definition 2.3} For a symplectic oriented cycle $\Si\subset S$,
 \roster
 \item The centre $c_\Si$ of the disc $D_\Si$ is called the {\it centre}
of $\Si$;
 \item the radius $r_\Si$ of the disc $D_\Si$ is called the {\it radius} of
$\Si$;
 \item $\Si$ is $\ep$-AG if
 $$
 c_\Si=0 \quad \text{and} \quad r_\Si<\ep.
 \tag2.54
 $$
 \endroster
\enddefinition

Donaldson proved the following analogue of Kodaira's embedding theorem (see
[D1]):

 \proclaim{Theorem 2.1} Let $L$ be a line bundle on an aK surface $S$ with
zero canonical class, and $a$ a Hermitian connection on $L$ whose
curvature form is the symplectic form on $S$:
 $$
 F_a=\frac{i}{2\pi} \om.
 \tag2.55
 $$
Then for any $k \gg 0$ there exists a section 
$s\in \Ga (L^k)$ such that
 \roster
 \item the zero set $s\1(0)=\Si$ is a smooth oriented symplectic $2$-cycle;

 \item $c_\Si=0$;

 \item $r_\Si<1/\sqrt k$.
 \endroster
In particular, there are many smooth symplectic cycles very close to AG
cycles.
 \endproclaim

Now consider an algebraic surface $S$ with canonical class $K_S>0$ or
$K_S<0$. We would like to deform slightly the algebraic geometry, by
considering a distinguished family of submanifolds, namely, the oriented
2-cycles. Moreover these cycles are determined by their first order
infinitesimal behavior, that is, the cycles defined by properties of
their Gauss lifts (2.29).

Recall that a $\Spin^{\C}$ structure is a choice of lift of the
projective bundle $\proj W^+$ to a vector bundle $W^+$. There are two
canonical choices of $\Spin^{\C}$ structure:
 $$
 W^+=\Oh_S \oplus \Oh_S(\pm K_S),
 \tag2.56
 $$
and we have to consider two cases: $K_S>0$ or $K_S<0$. 

\subheading{General type \rm($K_S>0$)} In this case, we consider the
lift
 $$
 W^+(K_S)=\Oh_S \oplus\Oh_S(K_S).
 $$
Recall that any section $s\in \Ga (W^+)$ is called a {\it spinor field}.
In particular, if we consider a nonvanishing section $s$ (or section
vanishing along a ``divisor'') then we get a section of the projective
bundle $\proj W^+$, that is, an almost complex structure $I_s$ compatible
with our K\"ahler metric $g$. 

In the K\"ahler (algebraic) case we have a finite dimensional family of
{\it holo\-morphic} sections
 $$
 s\in H^0(\Oh_S \oplus \Oh_S(K_S))=\C \oplus H^0(\Oh_S(K_S))
 \tag2.57
 $$
We write $p_g=\dim H^0(\Oh_S(K_S))$ for the geometric genus of $S$, that
is, the complex dimension of $H^{2,0}(S)$. The family of spinor fields
(4.23) defines a family of almost complex structures
 $$
 S^{2p_g}=\C^{p_g} \cup \{\infty\setminus\text{ point}\},
 \tag2.58
 $$
because the projectivisation 
 $$
 \proj(\C \oplus H^0(\Oh_S(K_S)))=\C^{p_g} \cup \proj H^0(\Oh_S(K_S)))
 \tag2.59
 $$
where $\C^{p_g}=\{(1, s)\}$ is the space of nonvanishing sections
and $\proj H^0(\Oh_S(K_S))=\vert K_S \vert$ is the complete canonical
linear system, points of which give the same almost complex structure,
namely, conjugate to the original complex structure.

Now the Levi--Civita connection gives the Hermitian structure on the
canonical bundle $L_K=\Oh_S(K_S)$ and, similarly, the Hermitian structure
on $\C^{2p_g} \cup \{\infty\setminus\text{point}\}$ gives the standard
metric on this sphere. Thus we can identify our sphere with the dual
sphere
 $$
S^{2p_g}=(S^{2p_g})^*
 $$
So this sphere contains the ``equator''
 $$
 S^{2p_g-1}_e=\{z \bigm| \Vert z \Vert=1\}\subset \C^{2p_g}.
 \tag2.60
 $$
Interpreting $\C^{2p_g} \cup \{\infty\setminus\text{point}\}$
as the space of sections gives us the embedding
 $$
 i_{\can} \colon \proj W^+\to S^{2p_g} \times S
 \tag2.61
 $$
where, as usual, $W^+=W^+(K_S)^*=\Oh_S \oplus \Oh_S(-K_S)$. The composite
of this embedding and the projection of the trivial bundle $S^{2p_g}
\times S$ to the fibre gives the map
 $$
 \pr\colon \proj W^+\to S^{2p_g}.
 \tag2.62
 $$
Now for any cycle $\Si\subset S$, the composite of the Gauss map (2.43),
the projection $\pr_+$ (2.44) and the projection (2.48) defines the
complex phase map
 $$
 m_\C={\pr}\circ {\pr_+} \circ G(i) \colon \Si\to S^{2p_g}.
 \tag2.63
 $$
Now in terms of this phase map, one can define the analogues of the
special Lagrangian cycles known in the Calabi--Yau case.

\definition{Definition 2.4} A cycle $\Si\subset S$ is called a sdAG cycle
if $m_\C(\Si)$ is a point, or equivalently
 $$
 \dd m_\C=0.
 \tag2.64
 $$
\enddefinition

The image $m_\C(\Si)\in S^{2p_g}$ defines an almost complex structure on
$S$ compatible with the K\"ahler metric for which $\Si$ is a {\it
pseudoholomorphic curve}. Hence, in particular, if $\Si$ is smooth, its
genus is
 $$
g(\Si)=\half([\Si]^2+[\Si] \cdot K_S)+1.
 \tag2.65
 $$
Now we can generalise ``Lagrangian'' properties of cycles:

\definition{Definition 2.5} 1) A cycle $\Si$ is called {\it weakly
Lagrangian} (wLag for short) if
 $$
 m_\C(\Si)\subset S^{2p_g-1}_e,
 \tag2.66
 $$
where $S^{3p_g-1}_e$ is the equator (2.60);

2) $\Si$ is called a spLag cycle, if it is a sdAG cycle and
 $$
 m_\C(\Si)\in S^{2p_g-1}_e.
 \tag2.67
 $$
\enddefinition

The equator divides the target sphere of the complex phase map into upper
and lower hemispheres:
 $$
S^{2p_g}=D^+ \cup S^{2p_g-1}_e \cup D^-,
 \tag2.68
 $$
and the entire catalogue of definitions (such as Definition 2.3),
properties and facts (such as Theorem 2.1) can be repeated in this new
set-up.

 \remark{Remark} If we start with an aK surface, we should remark
that instead of holomorphic spinor fields we should consider Dirac
harmonic spinor fields, that is, solutions of the differential equations
 $$
 D_a s=0,
 $$
where $D_a$ is the Dirac operator of the $\Spin^{\C}$ structure $-K_S$
coupled to the Levi--Civita connection on $W^+$, extended by the
connection $a$ on the determinant line bundle $L_{-K}$. More precisely,
the solutions of the Seiberg--Witten equations give the target sphere of
the complex phase map (2.60), and we can extend all our constructions to
aK case.
 \endremark

Finally, it is quite easy to see what to do if $K_S<0$: we change the
sign of the canonical system $K_S\to-K_S$, getting the sphere 
 $$
 S^{2h^0(\Oh_S(-K_S))}
 \tag2.69
 $$
as the target sphere of the complex phase map. After that, we can repeat
all our constructions and definitions. However, recall that 
 $$
 h^0(\Oh_S(K_S))=p_g=\half(b^+_2-1)
 \tag2.70
 $$
is a purely topological invariant of $S$, but $h^0(\Oh_S(-K_S))$ isn't a
topological invariant at all. Only for del Pezzo surfaces we can say that 
 $$
 h^0(\Oh_S(-K_S))\le10,
 \tag2.71
 $$
and
 $$
 h^0(\Oh_S(-K_S))=10-b_2^- .
 \tag2.72
 $$

\subheading{Remark} In the case 
 $$
h^0(\Oh_S(K_S))=h^0(\Oh_S(-K_S))=0,
 \tag2.73
 $$
we can't propose any good way of deforming algebraic geometry. This is
the rigid case.

 \head $\S$3. Slightly deformed Algebraic Geometry \endhead

Now we want to say why we need to deform ``algebraic geometry'' in this
style. But before starting on such explanations, we have to say what
algebraic geometry is, and what kind of questions we would like to discuss
in this set-up.

In algebraic geometry, we have a perfect theory of curves, a rather less
good theory of surfaces, an even worse theory of threefolds, and so on.
But now it is quite reasonable to consider the geometry of pairs
$C\subset S$, where $C$ is a smooth algebraic curve of genus $g$ and $S$ a
smooth algebraic regular surface containing $C$. Every irreducible
component $\sM$ of the moduli space of such pairs has a map
 $$
 f \colon \sM\to M_g \times \sM_s,
 \tag3.1
 $$
where $M_g$ is the moduli space of curves of genus $g$ and $\sM_s$ the
irreducible component of the moduli space of algebraic surfaces
containing $S$ as a point. Passing over the usual compactification
procedure, the fundamental class $[f(\sM)]$ gives cohomological
correspondences
 $$
 H_*(M_g)\to H^*(\sM_s) \quad \text{and} \quad H_*(\sM_s)\to H^*(M_g),
 \tag3.2
 $$
which are the ``topological'' part of our interest. On the other hand,
projections to components of the target space of $f$ define special
algebraic subvarieties of moduli spaces
 $$
 {\pr_g} \circ f (\sM)\subset M_g \quad \text{and} \quad {\pr_s} \circ f
 (\sM)\subset \sM_s,
 \tag3.3
 $$
the first of which is a very interesting subvariety of the moduli space of
curves of genus $g$ and the second usually coincides with the component
$\sM_s$, and the fibre of this projection
 $$
 ({\pr_s} \circ f)\1(S)=\vert C \vert
 \tag3.4
 $$
is the complete linear system of curves on $S$, the natural
compactification of the space of all deformations of $C$ inside $S$.

 \example{Example 1: \rm K3 surfaces} Let $S$ be a K3 surface and
suppose that the genus $g$ of the curve $C$ is $\ge12$. Then the
algebraic cohomology class $[C]$ defines a quasipolarisation of $S$, and
defines a component of moduli spaces $MK3_{[C]}$ of dimension 19:
 $$
 {\pr_s} \circ f (\sM)=MK3_{[C]}.
 \tag3.5
 $$
On the other hand, the compactified moduli space of deformations of $C$
inside a fixed K3 surface $S$ is the projective space
 $$
 \vert C \vert=\proj^g.
 \tag3.6
 $$
Thus
 $$
 \dim_\C \sM=g+19.
 \tag3.7
 $$
But, on the other hand, $({\pr_g}\circ f)(\sM)$ is a proper subvariety
of the moduli space $M_g$ of curves of genus $g$ and
 $$
 \codim_\C(({\pr_g} \circ f)(\sM))=2g-22
 \tag3.8
 $$
Thus the algebraic geometric problem is to describe this proper subvariety
in the moduli space of curves.
 \endexample

 \remark{Amazing observation} A generic algebraic curve of genus $>11$
lying on a K3 surface ``remembers'' the surface. Recall Mukai's recipe to
reconstruct the K3 surface $S$ containing a curve $C$ in terms of the
geometry of $C$ (for odd genus $g$).
 \endremark

 \proclaim{Mukai's recipe} 1) Consider the moduli space $\sSU_C(2, K_C)$
of semistable vector bundles on $C$ of rank 2 with $c_1$ the canonical
class. Inside this moduli space, consider the Brill--Noether locus
 $$
 \sSU_C(2,K_C,\half(g-1))=
 \bigl\{E\in\sSU_C(2,K_C)\bigm|h^0(E)\ge\half(g-1)\bigr\}. 
 \tag3.9
 $$
If $C$ is a general curve in $({\pr_g}\circ f)(\sM)$ then this locus is a
K3 surface $S'$.

2) The moduli space 
 $$
M_{S'}(2, \Theta \vert_{S'}, \half(g-5))=S
 \tag3.10
 $$
of rank $2$ torsion free sheaves with $c_1=\Theta \vert_{S'}$ (the
restriction of the theta divisor on the moduli space of vector bundles)
and $c_2=\half(g-5)$ is a K3 surface, namely the actual one containing $C$.
 \endproclaim

 \subheading{Remark} In this construction, Mukai realises his philosophy:
on the moduli space of vector bundles on a curve, the Brill--Noether level
plays the role of the second Chern class for vector bundles on a surface.

 \subheading{Basic classes} One has the same type of problem when
realising special integral 2-dimensional cohomology classes by special
geometric submanifolds, for example by (pseudo-)holomorphic curves. The
main example is following: the underlying smooth structure of an algebraic
surface with $p_g>0$ determines certain basic 2-dimensional cohomology
classes, the {\it Kronheimer--Mrowka} and {\it Seiberg--Witten} classes:
 $$
 \{\ka_i\}_{K-M}=\{\ka_i\}_{Z-W}\subset H^2(S,\Z).
 \tag3.11
 $$
These sets of classes are invariant under the diffeomorphism group and can
be realised as algebraic curves in every algebraic structure on $S$. The
set of these classes contains the canonical class $K_S$ realised by an
effective curve $C$. The genus of this curve is a purely topological
invariant of our surface:
 $$
 g(C)=2\chi+3 \si +1,
 \tag3.12
 $$
where $\chi$ is the Euler characteristic of $S$ and $\si$ its signature.
But the moduli of such curves is never generic. Indeed, the normal sheaf
of the canonical curve 
 $$
\Oh_C(C)=\Oh_C(\theta); \quad 2\theta=K_C, \quad h^0(\Oh_C(\theta))=p_g-1
 \tag3.13
 $$
is a theta characteristic on $C$ and if $p_g>2$, this theta characteristic
is special, that is, its theta constant vanishes. In this case one can
check that
 $$
\codim_\C ({\pr_g}\circ f) (\sM)=p_g-2
 \tag$3.13'$
 $$

On the other hand, the virtual dimension of the space of solutions of
Seiberg--Witten equations for a generic metric equals 0, but for a
K\"ahler metric, its actual geometric dimension equals $2(p_g-1)=b_2^+-2$.
Thus K\"ahler metrics are nongeneral and algebraic geometry isn't
``transversal'' for most problems of differential and symplectic
geomeries, for the theory of quantum multiplications and many others.

To describe the ``level of defectiveness'' of algebraic geometry, it is
enough to remark that the crucial point is the irregularity of the local
structure of the theory of deformations of curves inside a fixed
algebraic surface:

 \proclaim{Observation 3.1} On an algebraic surface $S$ with $p_g>0$, the
normal bundle $N_{C\subset S}=\Oh_C(C)$ of an algebraic curve $C$ is
almost always irregular, that is:
 $$
h^1(N_{C\subset S})>0.
 \tag3.14
 $$
 More precisely, $h^1(N_{C\subset S})=0$ only holds in two cases:

1) $C$ is a fixed component of the complete canonical linear system;

2) $C$ is a multiple fibre of an elliptic pencil.
 \endproclaim

 \remark{Remark} This property of normal bundles was observed by
Castelnuovo and Enriques 100 years ago as ``superabundance'': too many
curves with respect to points. In modern times it was remarked by
Donaldson in his seminal paper [D1].

The geometric meaning of this effect is quite simple: the virtual tangent
space to the space $\vert C \vert$ of deformations of curve $C$ in $S$ at
the point $C$ is
 $$
 T\vert C \vert_C=H^0 (N_{C\subset S}),
 \tag3.15
 $$
and the space of obstructions is
 $$
 H^1(N_{C\subset S}).
 \tag3.16
 $$
However, there are no genuine obstructions: every infinitesimal
deformation extends to a geometric deformation (just like under
deformations of Calabi--Yau manifolds): one has the exact sequence
 $$
H^0(\Oh_S(C))\to H^0(N_{C\subset S})\to H^1(\Oh_S)=0.
 \tag3.17
 $$
 \endremark

The final reason for wanting to slightly deform algebraic geometry is
purely arithmetic. Every smooth algebraic curve $C$ over $\Q$ is
traditionally considered as an algebraic surface fibred over $\Spec\Z$,
that is, as a pencil of curves. This pattern of thinking predicts to
consider rational points of $C$ as sections $s_1, \dots, s_N$ of this
pencil. By the Arakelov theorem squares of this sections are negative.
Suppose for a minute that these squares are $-1$. Then one can blow it
down to points $p_1,\dots,p_N$ on the surface $S$ with the curve $C$
(generic fibre) such that the normal bundle
 $$
 N_{C\subset S}=\Oh_C (p_1+\dots+p_N).
 \tag3.18
 $$
If $g(C)>1$ then $h^1(\Oh_C (p_1+\dots+p_N))>0$ and hence
 $$
 N<2g-2.
 \tag3.19
 $$
But this estimate is too good to be true. Thus our standard pattern is
wrong.

 \remark{Remark} The interesting fact is that this pattern was good
enough to prove such ``coarse'' fact as the Mordell conjecture but using
in Miyaoka's approach to prove Fermat's last theorem has shown that it
was wrong.
 \endremark

What we can do to avoid the speciality of moduli of curves which move in
algebraic surfaces, that is, the properness of the image $\pr_g$
projection (3.3)? The strong remedy is to consider all aK surfaces. In
this purely symplectic geometry set-up, one has to consider the space of
pairs $C\subset S$, where $S$ is any aK surface. Then the flexibility of
symplectic geometry gives immediately the result:

 \proclaim{Proposition 3.1} For aK surfaces, the image of the projection
$p_g$
 $$
({\pr_g} \circ \sM)=M_g
 $$
is the whole moduli space of curves of genus $g$.
 \endproclaim

Indeed, for any symplectic smooth oriented 2-cycles $\Si$ on
4-dimensional symplectic manifold $S$ and any complex structure on $\Si$
one can construct a compatible almost complex structure on $S$ such that
$\Si$ is a pseudoholomorphic curve for this almost complex structure
with induced complex structure (see, for example [M. Gromov ``Soft and 
hard symplectic geometry'', Proc.\ ICM, Berkeley, 1986, AMS, 1987,
81--98]). However we get a problem with the target space $\sM_s$ of the
second projection $\pr_s$: the ``moduli spaces'' of aK structures are
infinite dimensional. But fortunately, the question of describing fibres
of this second projections, that is, the ``complete linear systems'' of
pseudoholomorphic curves on aK surface is correct (because the
linearised operator is elliptic).

To be in the finite dimensional set-up, one has to fix some class of 
finite dimensional subspaces of $\sM$ of aK structures. And here we want
to change tack and to consider the slight deformations of almost complex
structures related to K\"ahler {\it metrics}.

We would like to deform algebraic geometry slightly in such a way that
the theory of curves will be preserved completely and in a pair $(C\subset
S)$ we have to consider, instead of an algebraic curve $C$, a smooth sdAG
cycle $\Si$ of fixed phase
 $$
p=m_\C(\Si)\in S^{2p_g}.
 \tag3.20
 $$
Every ``irreducible'' component of the moduli space $\sM$ of such pairs
$(\Si\subset S)$ defines a map $f$ (3.1), a cohomological correspondence
(3.2) (because every sdAG cycle admits a complex structure) and a subspace
(3.3) in the moduli space of curves of genus $g$. In particular one has to
describe the global structure of the complete space $\vert \Si \vert$ of
deformations of sdAG cycle $\Si$ in $S$. Such type cycles are
pseudoholomorphic curves for almost complex structure $S^p$ given by the
point $p$ (3.20). Of course, our initial complex structure is $S^0$ and
its complex conjugate is $S^{\infty}$. The local deformation theory of sdAG
cycles is quite good, like the theory of local deformations of spLag cycles
(see [H-L]). 

Now fix an algebraic surface $S$ and consider a smooth sdAG cycle $\Si$
of the phase $p$ (3.20). Let 
 $$
{\Si}\in \vert [\Si] \vert^p
 \tag3.21
 $$
be the moduli space of deformations of $\Si$ as a sdAG cycle. Then there
exists a deformation $\de$-complex such that the tangent space is
given by
 $$
 T\vert [\Si] \vert^p_\Si
 =H^0_{\de} (N_{\Si\subset S})=H^0(\Si, N_{\Si\subset S}).
 \tag3.22
 $$
(Recall that $\Si$ admits a complex structure, in which the normal bundle
$N_{\Si\subset S}$ admits a holomorphic structure. The last space is a
coherent cohomology group of a holomorphic vector bundle.)

In the same vein, the space of obstructions of infinitesimal deformations
of $\Si$ as a sdAG cycle is
 $$
 H^1_{\de}(N_{\Si\subset S})=H^1(\Si,N_{\Si\subset S}),
 \tag3.23
 $$
where the last space is again the space of coherent cohomology of the
holomorphic vector bundle $N_{\Si\subset S}$ on the complex curve $\Si$.

There are two transversality results. The first one was proposed by
Donaldson in [D1]:

 \proclaim{Theorem 3.1} 1) For a generic aK structure, the obstruction
space to local deformations of any pseudoholomorphic curve vanishes
 $$
H^1_{\de} (N_{\Si\subset S})=H^1(\Si, N_{\Si\subset S})=0.
 \tag3.24
 $$

2) The dimension of the local deformations space is given by
 $$
 \dim T\vert[\Si]\vert^p_\Si=H^0_{\de}(N_{\Si\subset S})
 =H^0(\Si,N_{\Si\subset S})=2(\Si^2+1-g).
 \tag3.25
 $$
 \endproclaim

We can add to this the same type ``transversality result'':

 \proclaim{Theorem 3.2} If $S$ is a surface with $p_g>0$ then for
generic point $p\in S^{2p_g}$ of the target space the complex phase map
and any sdAG cycle of the phase $p$ we have equalities (3.24) and (3.25).
 \endproclaim

The idea of the proof is contained in the same paper of Donaldson [D1]. 
First of all we would like
to consider the whole family 
 $$
 \Vert \Si \Vert=\bigcup_{p\in S^{2p_g}} \vert [\Si] \vert^p
 \tag3.26
 $$
of sdAG cycles for all $p\in S^{2p_g}$ and to prove that this family is
smooth and of the right dimension (that is, the virtual dimension). After
using standard arguments we get the statement of the theorem. 

The whole family of sdAG cycles doesn't fibre over the target of the phase
map space because one sdAG cycle can be sdAG in more than one complex
structures. But there exists a correspondence, or a space of pairs
 $$
\sM=\{(\Si, p) \bigm| \Si \text{ is sdAG with $m_\C(\Si)=p$}\}\subset
\Vert \Si \Vert \times S^{2p_g},
 \tag3.27
 $$
with two projections
 $$
 \sM @> \pr_t >> S^{2p_g} \text { and } \sM @> \pr_c>> \Vert \Si \Vert 
 \tag3.28
 $$
and fibres of the second projection are 
 $$
 \pr_t\1(p)=\vert [\Si] \vert^p.
 $$

\subheading{Example. The canonical system} Consider the canonical class
$K_S$ of $S$, and realise it as sdAG cycles for all $p\in S^{2p_g}$. Then
one has

 \proclaim{Proposition 3.2}
 $$
 \Vert K_S \Vert=\vert K_S \vert^0.
 \tag3.27
 $$
In particular, the image of the map ${\pr_g} \circ f (\sM)$ (3.3) is the
same as for canonical curves of the algebraic surface $S^0$ (see (3.13)).
 \endproclaim

To describe the space of pairs (3.27) we must blowup the point $(1,0)$ in
the projective space $\proj(\C \oplus H^0(\Oh_{S^0} (K_S))$ (2.43),
(2.45):
 $$
\si \colon \proj(\Oh \oplus \Oh(H))\to \proj(\C \oplus H^0(\Oh_{S^0} (K_S)
 \tag3.28
 $$
and the second projection of this blown up projective space 
 $$
\pi \colon \proj(\Oh \oplus \Oh(-H))\to \vert K_S \vert^0,
 \tag3.29
 $$
where $\Oh$ is the structure sheaf of the projective space 
 $$
 \proj H^0(\Oh_{S^0} (K_S))=\vert K_S \vert^0
 \tag3.30
 $$
and $\Oh(-H)$ is the Hopf line bundle on this projective space.

Then it is easy to see that
1) the space of pairs (3.27) is 
 $$
 \sM=\proj(\Oh \oplus \Oh(-H)),
 \tag3.31
 $$
and the map to the target sphere of the complex phase map 
 $$
 \pr_t \colon \sM=\proj(\Oh \oplus \Oh(H))\to S^{2p_g}
 \tag3.32
 $$
is the blowdown of the sections $\proj(\Oh)$ and $\proj(\Oh(-H))$ in the
projective bundle (3.31) over the projective space $\vert K_S \vert^0$:
 $$
 \pr_t(\proj\Oh)=0\in S^{2p_g} \quad \text{and}\quad
 \pr_t(\Oh(-H))=\infty\in S^{2p_g}.
 \tag3.33
 $$
{From} this picture one get immediately

 \proclaim{Proposition 3.3} 1) For any $p\in S^{2p_g} \setminus (0
\cup\infty)$ on aK surface $S^p$, the canonical class admits a unique
representation as pseudoholomorphic curve 
 $$
 \vert K_S \vert^p=\text{a point};
 \tag3.34
 $$

2) For every curve $C\in \vert K_S \vert^0$ there exists a 2-sphere 
 $$
 S^2=\pr_t (\pi\1(C))\in S^{2p_g}
 \tag3.35
 $$
of almost complex structures in every of which $C$ is the
pseudoholomorphic realisation of the canonical class.
 \endproclaim 

\head $\S$4. The complex 3-dimensional case \endhead

Now we would like to extend the constructions we have described to the
case an aK threefold $X$ given by an Hermitian triple $(\om,I,g)$. Here we
describe simple facts and constructions from linear algebra of such
structure over one point. Afterwards, in the next section, we globalise
these constructions to the geometry of the tangent bundle and special
subspaces of spaces of tangent directions to define new collection of
submanifolds and hence, new collections of geometries of complex
threefolds.

In the even dimensional case, one has two volume forms as the formula
(2.9). In the odd dimensional case, one has to use the special trick to
reduce 3 to 2.

Over any point $p\in X$ we have as the fibre of the tangent bundle
$TX_p=\C^3$ with the standard Hermitian triple $(\om_p, I_p, g_p)$ with
the constant symplectic form $\Span{\ ,\,}=\om_p$, the constant Euclidean
metric $g_p$, giving the Hermitian triple $(\om_p,I_p,g_p)$ and the
oriented Lagrangian Grassmannian is a subspace of the oriented Grassmannian
 $$
(\LaGr)_p\subset\orGr(3, TX_p)
 \tag4.1
 $$
of codimension 3. Moreover, instead of the complex Grassmannian (2.13) one
has the subspace
 $$
 \{V\in\orGr(3, TX_p) \bigm| V \text{ contains a complex direction }
z\in \proj TX_p\}
 \tag4.2
 $$
of subspace of 3-dimensional planes containing complex directions. Just as
in the 2-dimensional case, for every $V\in\orGr(3, TX_p)$ there are two
possibilities:
 $$
\dim\Span{V,I(V)}=6.
 \tag4.3
 $$
This is the general case. In particular every Lagrangian $V$ has this
property. The second case is when
 $$
\dim\Span{V,I(V)}=4.
 \tag4.4
 $$
In this case $V$ containes a complex direction (see (4.2)) and defines the
flag
 $$
z_V\subset V\subset Z_V
 \tag4.5
 $$
where $z$ is the unique complex direction in $V$ and $Z_V=\Span{V,I(V)}$
(4.4).

We can distinguish two components of the set (4.2) by the orientation: the
orientation of $V$ may or may not be compatible with the complex
orientation of $z_V$. Now we can realise these components as level
submanifolds of a map of $\orGr(3,TX_p)$ to $S^2$ (just as in the
2-dimensional case) for complex and anticomplex directions. The following
constructions work for any dimensional case. We described its in \S2 (see
(2.18--2.34)) but for simplicity, we repeat it for the 3-dimensional case
again.

\definition{Definition 4.1} Any nonzero element of $w\in\bigwedge^3 \C^3$
is called a {\it complex orientation} of $\C^3$.
 \enddefinition

Let $w\in\bigwedge^3 \C^3$ be a fixed complex orientation of $\C^3=TX_p$.
Let us fix some orthogonal basis $e_1, e_2, e_3, I(e_1), I(e_2), I(e_3)$ of
$TX_p$ such that
 $$
 w=(e_1+I(e_1))\wedge (e_2+I(e_2))\wedge (e_3+I(e_3)).
 \tag4.6
 $$ 
If we fix the volume form $w$ such a way that $|w|=1$ then we have to fix
a phase $e^{i\fie}$ only.

Now in any oriented subspace $V\in\orGr(3, TX_p)$, consider an orthogonal
basis $v_1, v_2, v_3$ and let $M_V$ be the linear transformation sending
$e_i$ to $v_i$ and $I(e_i)$ to $I(v_i)$. Again it is easy to see the
following:

 \proclaim{Proposition 4.1} 1) $M_s$ is a $\C$-linear map $TX_p\to TX_p$,
that is, a complex matrix in our basis;

2) $\det M_V=0$ iff $V$ contains a complex direction $z$. Of course $z_V$
is the kernel of $M_V$;

3) $V$ is Lagrangian iff $M_V\in\U(3)$. 
 \endproclaim

Again we have the map
 $$
d \colon\orGr(3, TX_p)\to \C
 \tag4.7
 $$
sending $V$ to the determinant of $M_V$ with the same inequality (2.25),
but the one difference is in (2.26): one has, of course, the following:

 \proclaim{Hadamard's Lemma}
 $$
V\in(\LaGr)_p \iff \vert d(V) \vert=1.
 $$
(See for example [H-L], Lemma 1.9.)
 \endproclaim
Thus the oriented Lagrangian Grassmannian is
 $$
(\LaGr)_p=d\1(S^1)
 \tag4.8
 $$
However, the dimension of all others fibres is 7, whereas the dimension of
$d\1(e^{i\fie})$ is 5 -- an extremely small number!

Now we can finish our construction. Again we have the decomposition
 $$
 \orGr(3, TX_p)=\orGr(3, TX_p)^+ \cup d\1(\{\vert d(V) \vert=1\})
 \cup\orGr(3, TX_p)^-
 \tag4.9
 $$
with the inversion orientations map (2.31). Recall again that the middle term is
 $\LaGr$.

Then the glueing of two maps $d$ on $\orGr(3, TX_p)^+$ and $-(\overline{d})\1$ which
are equal on $d\1(\{\vert d(V) \vert=1\})$ gives the map
 $$
\det_{I} \colon\orGr(3, TX_p)\to \proj^1= S^2_w
 \tag4.10
 $$
(just compare with (2.30--2.31)). But the difference between the
2-dimensional and 3-dimensional case is that this map {\it is not smooth
along} $(\LaGr)_p$!

To make it smooth we have to ``blowup'' $(\LaGr)_p$ inside $\orGr(3,
TX_p)$. We do it after the couple of remarks about complex orientations.

Of course the map $\det_{I}$ depends on $I$ and on a complex orientation
$w$ but the fibres
 $$
(\det_{I})\1(0) \quad \text{and} \quad (\det_{I})\1(\infty)
 \tag4.11
 $$
do not depend on a complex orientation. Only the identification of
 the target space $\proj^1=S^2_w$ depends
on $w$.

Indeed the space of complex orientations 
 $$
 \{w\}=\C^* \quad \text{ (or $e^{i\fie}$ if } |w|=1)
 \tag4.12
 $$
acts on the complex sphere $S^2=\C \cup\infty$ by the usual multiplication: 
for any $\ep\in \C^*$
 $$
 \det_{I, \ep \cdot w}=\ep \cdot \det_{I, w}
 \tag4.13
 $$
and one get
 \proclaim{Proposition 4.2} 1) All fibres of the map $\det_I$ (4.10) over 
$S^2\setminus(0 \cup\infty \cup S^1)$ are diffeomorphic:
 $$
 (\det_{I})\1(z)=(\det_{I})\1(z') \quad\text{for }
 z, z' \ne0,\infty, e^{2\pi i\fie}.
 \tag4.14
 $$
 \endproclaim

Of course one can send the target space $S^2$ of this map to the interval
$[0,\pi]$ by the map $h$ (2.12) such a way that one get the ``K\"ahler
angle'' map:
 $$
 \al_{I, w}=h \cdot \det_{I, w} \colon\orGr(3, TX_p)\to [0,\pi].
 \tag4.15
 $$
(just as in the 2-dimensional case).

Moreover, as in the 2-dimensional case,
 $$
 \dim\al_I\1(z)=8, \quad \dim\al_I\1(e^{i\fie})=5 \quad
 \text{and} \quad \dim\al_I\1(0)=7,
 \tag4.16
 $$
and for $\al_I\1(z)$ (4.14) there exists an extra phase map
 $$
 \al_I\1(z)\to S^1,
 \tag4.17
 $$
whose restriction to the oriented Lagrangian Grassmannian is equal to the
standard phase map
 $$
 \det \colon\LaGr(3, TX_p)\to S^1_p
 \tag4.18
 $$
(see the begining of our story in \S1).

Every $V\in\orGr(3, TX_p)$ can contain one complex direction
only, thus
 $$
 \rk M_V \ge 2.
 \tag4.19
 $$
 \remark{Remark} Hence instead of $M_V$ we can consider $\adj M_V$ and the
image
 $$
 \im (\adj M_V)=z_V; \quad \ker (\adj M_V)=Z_V
 \tag4.20
 $$
 \endremark

So if $\det M_V=0$ then one has the orthogonal decomposition
 $$
 V=\ker M_V \oplus (\ker M_V)^{\perp_V}
 \tag4.21
 $$
where $\perp_W$ is the orthogonal to subspace of a Euclidien space $W$.
Thus $(\ker M_V)^{\perp_V}$ is an oriented 1-subspace in $V$. But the plane 
 $$
 \Span{(\ker M_V)^{\perp_V}, I((\ker M_V)^{\perp_V})}
 \tag4.22
 $$
is complex. This is a description of the 
canonical flag $(z_V\subset V\subset Z_V)$ (4.5).

So one get the map
 $$
f \colon \al_I\1(0)\to F_{1,2}^{\C} (\proj (TX_p))=\proj T \proj(TX_p)
 \tag4.23
 $$
to the complex flags of type (1,2) in $TX_p=\C^3$. Obviously a fibre of
this map is $S^1$. More precisely, the $\proj^1$-bundle which is the target
space of the map $f$ (4.23) has the tautological bundle $H$ (the
Grothendick line bundle of a projectivisation) equipped with the Hermitian
metric. The unit circle bundle of this line bundle is the source of (4.23):
 $$
\al_I\1(0)=S^1(H).
 \tag4.24
 $$
Recall that this is a description of any fibre of the map $\det_I$ expect
for $\{e^{i\fie}\}=S^1$.

Now if one consider the projection to $\proj(TX_p)$ then the fibre of this
map over a point $z$ is the space of real rays in
the vector space $\C^2=\C^3 / z$. Thus this fibre is the unit 3-sphere
 $$
\pr_\C\1(z)=S^3\subset \C^3 / z.
 \tag4.25
 $$
Using the Euler exact sequence it is easy to see that full fibre of $\al_I$ over $0$ is
 $$
\al_I\1(0)=S^3 (T \proj (TX_p) (-1))
 \tag4.26
 $$
that is, the unit 3-spheres fibration of the twisted tangent bundle
$T\proj (TX_p) (-1)$ of the complex projectivisation of the tangent space
to $X$ at $p\in X$. Of course
 $$
 S^3 (T \proj (TX_p) (-1))=S^1(H)
 \tag4.27
 $$
and fibrewise 
 $$
S^3\to \proj^1
 \tag4.28
 $$
is the Hopf fibration.

 \remark{Remark} This is the first difference between the complex dimension
2 case and the 3-dimension case: for 2-dimension case the fibre of $\al$
over $0$ is just the projectivisation of the tangent space $\proj W_-$ at
$p$ (see (2.13)).
 \endremark

The second is that the preimage 
 $$
 \al_{I, w}\1(\pi/2)=(\LaGr)_p 
 \tag4.29
 $$
has extremely small dimensional.

\definition{Definition 4.2} 1) An oriented 3-dimensional subspace $V\in\orGr(3, TX_p)$
 is called 3/2-{\it pseudoholomorphic} (3/2ps for short) if 
 $$
 \al_I (V)=0;
 \tag4.30
 $$

2) it is called 3/2-{\it anti-pseudoholomorphic} (3/2aps for short) if
 $$
 \al_I (V)=\infty.
 \tag4.31
 $$

3) a pair $(v\subset V)$ where $v$ is oriented 1-dimensional subspace of
$V$ is called {\it super Lagrangian} (superLag for short) if
 $$
 V\in (\LaGr)_p.
 \tag4.32
 $$
and super spLag if $V$ is spLag.
\enddefinition
 \remark{Remark} 1) We know that the fibres 1) and 2) don't depend on a
complex orientation $w$, but 3) depends on $w$.

2) For the definition of $\LaGr$ we have to use the form $\om$ which
depends on $I$. So it will be quite correct to write $\LaGr(I)$.
 \endremark 

Now we can describe the map $\al_I$ (4.7) in other way using the form $\om$ which is defined by
I (and our metric). Let us consider the
open set
 $$
 \orGr(3,T X_p)_{\ne0}=\orGr(3, T X_p)-(\LaGr)_p.
 \tag4.33
 $$
defined by
 $$
 V\in \orGr(3, T X_p)_{\ne0} \iff \om \vert_V \ne0.
 $$
Then any such $V$ admits the decomposition
 $$
V=\ker \om \vert_V \oplus (\ker \om \vert_V)^{\perp_V}
 \tag4.34
 $$
and we {\it fix the orientation of the plane}
 $$
t_V=(\ker \om \vert_V)^{\perp_V}
 \tag4.35
 $$
{\it such
a way that the volume} $\om \vert_t$ {\it is positive}.

Then we get the map
 $$
 \be_{\om} \colon \orGr(3, T X_p)_{\ne0}\to \orGr(2, T X_p)
 \tag4.36
 $$
sending $V$ to $t_V$ (4.25).

But for $\orGr(2, T X_p)$ we have the {\it classical} K\"ahler angle
map
 $$
\al_p \colon \orGr(2, T X_p)\to [0,\pi]
 \tag4.37
 $$
sending a plane $t$ to the same K\"ahler angle (2.9).

 \remark{Remark} This construction only depends on the {\it conformal class}
of $\om$.
 \endremark

What we have to do now it is just ``blowup the oriented Lagrangian
Grassmannian $(\LaGr)$ inside the Grassmannian'': consider the space of pairs
 $$
 \wave{\orGr(3, T X_p)}=\{(v\in V) \bigm| v \text{ is the kernel of } \om
 \vert_V\}
 \tag4.38
 $$
and its projection to $V$:
 $$
\si \colon \wave{\orGr(3, T X_p)}\to \orGr(3, T X_p).
 \tag4.39
 $$

 \proclaim{Proposition 4.3} \roster
 \item
 The map $\si$ is an isomorphism over $\orGr(3, T X_p)_{\ne0}$.

 \item
 $$
 \si\1(\LaGr)=S^2(U)
 \tag4.40
 $$
is the unit sphere bundle of the universal bundle $U$ over the oriented
Lagrang\-ian Grassmannian $\LaGr$. This space is the moduli space of super
Lagrangian $3$-directions (4.32).

 \item The map (4.36) has a smooth extension to
 $$
 \be_{\om} \colon \wave{\orGr(3, T X_p)}\to \orGr(2, T X_p);
 \tag4.41
 $$
and for any Lagrangian super direction $v\in V$
 $$
 \gather
 \be_{\om}(v\subset V)=v^{\perp_V}, \quad\text{that is,}\quad
 t_{(v\subset V)}=v^{\perp_V} \\
 \text{and}\quad V \text{ is } 3/2-ph \implies \be_{\om}(V)=z_V.
 \endgather
 $$

 \item
 $$
 \be_{\om}\si\1(\LaGr)=\LaGr(2,TX_p)\subset \orGr(2,TX_p),
 $$
and any fibre $\be_{\om}\1(t)\subset S^3(Q_4)_t$, where $Q_4$ is the
universal quotient bundle on $\orGr(2, T X_p)$; let 
 $$
 \proj TX_p\subset \orGr(2, T X_p)
 $$
be the projective plane of complex directions in the tangent space then
 $$
 \be_{\om}\1(t)=I(t) \cap S^3(Q_4)_t=S^1 \iff
 t\in \orGr(2, T X_p)- \proj TX_p
 $$
and
 $$
 \be_{\om}\1(t)=S^3(Q_4)_t \iff t\in \proj TX_p.
 $$

 \item The composite of the maps (4.37) and (4.41) gives a smooth map
 $$
 \wave{\al_I} \colon \wave{\orGr(3, T X_p)}\to [0,\pi].
 \tag4.42
 $$
Its geometric meaning is quite simple: $\wave{\al_I}$ {\it sends $V$ to
the K\"ahler angle of} $t_V$. In particular $\wave{\al_I}((v\in V))=\pi/2$.
 \endroster
 \endproclaim

Our map depends on the complex structure $I$. Now, {\it what happens if we
change the initial complex structure slightly?}

In our odd dimensional case, an oriented 3-subspace $V\in\orGr(3,TX_p)$
doesn't determine a new complex structure on $TX_p$ compatible with the
metric $g$. But a flag $(t\subset V\subset T)$ where $t$ is an oriented
2-plane in $V$ and $T$ is an oriented 4-subspace determines a new complex
structure on $TX_p$ compatible with the metric $g$. To see this, consider
the orthogonal decomposition into 2-dimensional subspaces
 $$
TX_p=t \oplus t^{\perp_T}\oplus T^{\perp}
 \tag4.43
 $$
and put the standard complex structure on each 2-subspace of this
decomposition. Then one get a new complex structure 
 $$
I_{(t\subset V\subset W)}
 \tag4.44
 $$
on $TX_p$ such that $t$ is a complex direction and $T$ a complex
subspace. 

Now the description of the space of all flags $(t\subset V\subset T)$ is:
 $$
F_{2,3,4}=S^2(U) \times_{\orGr(3, TX_p)} S^2 (Q_3)\to\orGr(3, TX_p).
 \tag4.45
 $$
where $U$ is the tautological bundle on the Grassmannian, $Q$ the universal
quotient bundle and $S^2(*)$ the unit sphere bundle of $*$.

Every new complex structure (4.44) determines a map
 $$
\al_{I_{(t\subset V\subset T)}} \colon\orGr(3, TX_p)\to [0,\pi] 
 \tag4.46
 $$
which isn't smooth along the oriented Lagrangian Grassmannian
$\LaGr(I_{(t\subset V\subset T)})$, and a map $\wave{\al}_{I_{(t\subset
V\subset T)}}$ (4.42) which is the regularisation of $\al_I$. 

Again one has the subspace of 3/2-pseudoholomorphic oriented 3-subspaces
 $$
\al_{I_{(t\subset V\subset W)}}\1(0)\subset\orGr(3,TX_p).
 \tag4.47
 $$
and so on. 

\head \S5. Complex structures and globalisations \endhead

It is now time to describe the space of complex structures on $TX_p$
compatible with our metric $g_p$. To do this, consider the
complexification of the metric quadric $g_p^{\C}$ in
$TX_p^{\C}=TX_p\tensor\C$; any compatible complex structure on $TX_p$ is
given by a maximal isotropic subspace $T^{1,0}\subset TX_p^{\C}$ with
respect to our quadric $g_p^{\C}$. For algebraic geometers it is quite
convenient to projectivise all geometric objects. So we have
 $$
G_p\subset \proj^5=\proj TX_p^{\C}
 \tag5.1
 $$
where $G_p$ is the smooth 4-dimensional quadric in $\proj^5=\proj
TX_p^{\C}$ of isotropic lines. The projectivisation of $T^{1,0}\subset
TX_p^{\C}$ is a projective plane in $G_p$. There are two systems of planes
in any nonsingular quadrics which one can interpret as the Grassmannian of
lines in $\proj^3_+$. Then one system of planes on $G_p$ is given by points
of this $\proj^3_+$ (as the set of lines through a point) and other one is
given by planes (as the set of lines in fixed plane) that is, points of the
dual space $(\proj^3_+)^*=\proj^3_-$.

 \proclaim{Proposition 5.1} These systems of planes on $G_p$ are
distinguished by the orientation of $X$.
 \endproclaim

Each of these spaces is the projectivisation of a spinor space at a point
$p\in X$:
 $$
\proj^3_{\pm}=\proj W^{\pm}; \quad \proj W^+=(\proj W^-).
 \tag5.2
 $$
The projective space $\proj TX_p^{\C}$ has a real structure, with respect
to which the metric quadric $G_p$ (4.34) is real:
 $$
\theta \colon G_p\to G_p
 \tag5.3
 $$
without fixed points, that is, without real points. Therefore $\theta$ must
send one system of planes on $G_p$ to the other:
 $$
 \theta \colon \proj^3_+\to \proj^3_-=(\proj^3_+)^*.
 \tag5.4
 $$
Thus $\theta$ determines a projective Hermitian structure on $\proj^3_+$.
One gets

 \proclaim{Proposition 5.2} The space of complex structures on $TX_p$ 
compatible with the metric $g_p$ and the orientation of $X$ is\/
$\proj^3=\proj W^+_p$. (Just as in the 2-dimensional case).
 \endproclaim

Now one has the map
 $$
c_p\colon F_{2,3,4}=S^2(U) \times_{\orGr(3, TX_p)} S^2 (Q)\to \proj W^+_p.
 \tag5.5
 $$
sending a flag $(t\subset V\subset W)$ to the complex structure
$I_{(t\subset V\subset W)}$ (compare with (2.6)). Here and over there $U$
and $Q$ are the universal subbundle and quotient bundle of the
Grassmannian.

 \proclaim{Proposition 5.3} The fibre
 $$
c_p\1(I_{(t\subset V\subset W)})=\al_{I_{(t\subset V\subset W)}}\1(0). 
 \tag5.6
 $$
 \endproclaim
(See (4.46--4.47)).

Our complex structure $I$ is the point of $\proj W^+_p$
 $$
I_p\in \proj W^+_p
 \tag5.7
 $$
and the conjugated complex structure 
 $$
\overline{I}=\theta (I)\in \proj W^-_p .
 \tag5.8
 $$

 But this point determines the projective plane
 $$
(\proj^2_{\overline{I}})_p\subset \proj W^+_p.
 \tag5.9
 $$

 \proclaim{Proposition 5.4} 1) Every pair $(I,I')$ of complex structures
on $TX_p$ has unique common complex direction.

2) The plane $\proj^2_{\overline{I}}$ is canonically isomorphic to the
projectivisation $\proj TX_p$ of the tangent space (in the complex
structure $I$).

3) Every oriented plane $t$ in $TX_p$ defines a line $l_t\in \proj
W^+_p$, that is, a point of the complex quadric $G_p$ (5.1). This map
 $$
s_p \colon\orGr(2, TX_p)\to G_p
 \tag5.10
 $$
is an isomorphism.
 \endproclaim

Indeed, every pair of points $I,I'\in \proj W^+_p$ determine two
projective planes 
 $$
\proj T^{1,0}_I, \proj T^{1,0}_{I'}\subset G_p\subset \proj TX^{\C}_p
 \tag5.11
 $$
(see (5.1)). These planes in the quadric have one intersection point
 $$
\proj T^{1,0}_I \cap \proj T^{1,0}_{I'}=z (I, I')\subset \proj TX^{\C}_p
 \tag5.12
 $$
which gives the common holomorphic direction. This gives 1).

To get 2), consider the map
 $$
z(I,\ ) \colon \proj^2_{\overline{I}}\to \proj T^{1,0}_I 
 \tag5.13
 $$
sending a plane $I'\in \proj^2_{\overline{I}}$ to the intersection point
$z (I, I') $. It is easy to see that this map is an isomorphism.

So one get the following interpretation of the right hand spinor space (5.2)
 $$
W^+_p=\C \cdot I \oplus TX_p
 \tag5.14
 $$
because of the equality $TX_p=T^{1,0}_I $.

Now for any point $p\in \proj^2_{\overline{I}}$ one has the projective
line 
 $$
 \Span{I,p}\subset \proj W^+
 \tag5.15
 $$
of complex structures.

 \proclaim{Proposition 5.5} 1) Every complex structure of 
the family of complex structures $\Span{I,p}$ has the same complex
direction common with $I$.

2) For every complex direction $p$ of a complex structure $I$ the family
$F_I(p)$ of complex structures every of which containes $p$ as a complex
direction is the family (5.15):
 $$
 F_I(p)=\Span{I,p}\subset \proj W^+
 \tag5.16
 $$
 \endproclaim

Now we can identify the family of complex structures $F_I(p)$ (5.16) with
$p$ given by the complex direction $t\subset T^{1,0}_I$ with the space of
complex structures on 
 $$
 \C^2=\R^4=T^{1,0}_I/t.
 \tag5.17
 $$
This can be identified with the target space of the map $\pr_+$ (2.6) for
the complex 2-dimensional case (see \S2). That is, we have the
identification
 $$
 \Span{I,p}=\proj^1_+
 $$
for $TX_p=\R^4$ and the map
 $$
 \pr_+ \colon\orGr(2,T\proj^2_{\theta(I)})\to \Span{I,p}
 \tag5.18
 $$

Now we have to switch on our form $\om$ to get the K\"ahler angle map
(2.11)
 $$
 \orGr(2, T \proj^2_\infty) @>\pr_+>> {\Span{I,p}} @>h>> [0,\pi];
 \tag5.20
 $$
see (2.7--2.13).

So under the identification of a point $p$ of the plane
$\proj^2_{\theta(I)}$ to the complex direction $t\subset T^{1,0}_I$ we get
an identification of the family $F_I(p)$ (5.15) with the family of complex
structure $\Span{I,p}=S^2$ on $\C^2=T^{1,0}_I/t$ (5.16).

Thus all the maps (2.12)
 $$
 h \colon\Span{I,p}\to [0,\pi]
 \tag5.21
 $$
are restrictions to lines $\Span{I,p}\subset\proj^3$ of the map
 $$
 \proj^3 @>\overline{\si}>> S^6_I @>h>> [0,\pi],
 \tag5.22
 $$
where $\overline{\si}$ is the blowdown of the projective plane
$\proj^2_{\theta}$ to the point $\infty$. Let us denote this composite by
 $$
h_I \colon S^6_I\to [0,\pi]
 \tag5.23
 $$
and call it the {\it height} function (for the 3-dimensional case).

This 6-sphere $S^6$ is a fibre of the Thom space of the complex vector
bundle $TX^I$ for the initial aK structure on $X$.

So in the 2-dimensional case we have the map (2.12) $\proj^1_+=S^2 @>h>> [0,\pi]$ and in
3-dimensional case we have the map (5.23) of a fibre of the Thom space
 $$
Th(TX^I_p) @>h_i>> [0,\pi].
 \tag5.24
 $$

Again we have two points $I=0$ and 
 $$
 \overline{I}=\theta(I)=\infty
 \tag5.25
 $$
such that all its inverse images are points:
 $$
h_I\1(0)=I ; \quad h_I\1(\infty)=\overline{I}
 \tag5.26
 $$
and all other fibre we can lift by $\overline{\si}\1$
(5.22) to $\proj^3$ and project (from $I$) to the 
projective plane $\proj^2_{\theta}$:
 $$
\psi=\overline{\si}\1\cdot \pr_I \colon h_i\1(\fie)\to \proj^2_{\theta}.
 \tag5.27
 $$
It is easy to see that this map is surjective and the fibre over an
$I$-complex direction $z\in \proj TX^I_p$
 $$
\psi\1(z)=h\1(\fie)
 \tag5.28
 $$
is the actual phase circle of the standard K\"ahler angle map (2.12)
(2-dimensional case) 
 $$
h \colon\Span{I,z}\to [0,\pi]
 \tag5.29
 $$
sending a complex structure $I'\in\Span{I,z}$ to the complex structure on
$\R^4=TX/z$ with K\"ahler angle $\fie$.

In particular $h\1(\pi/2)$ is the circle of complex structure
$\{I'\}\subset\Span{I,z}$ such that any $I'$-complex direction in
$\R^4=TX/z$ is Lagrangian (more precisely $I$-Lagrangian.

Now in the direct product 
 $$
 \orGr(2, TX_p) \times \proj^3_+
 \tag5.30
 $$
 one has subspace
 $$
 \bigwedge\suniv=\{(t, I) \bigm| \om_I \vert_t=0\}.
 \tag5.31
 $$
In particular for every $t\in\orGr(2, TX_p)$ we have
 $$
 \sI_t=(t, \proj^3_+) \cap\bigwedge\suniv\in \proj^3_+,
 \tag5.31
 $$
that is, the space of all complex structure for which $t$ is Lagrangian.

In the same vein, we have subspaces of $\proj^3_+$:
 $$
 \sI_V=\{I \bigm| \om_I \vert_V=0\}
 \tag5.32
 $$
and so on. We leave the problem of describing these subspaces as exercises
for the reader.

\subheading{Globalisation} In the 3-dimensional case, our pointwise
constructions have ``moduli''. Before globalising them, let us bring
together all observations coming from our local investigations.

 \roster
 \item First of all, the dependence of the map $\det_{l, w}$ on $I$ and on
the complex orientation $w$ (see (4.6)) was described by the formula
(4.13).

 \item The blowup of the Grassmannian $\wave{\orGr(3, TX_p)}$ (4.39) has
the regular map $\wave{\al_I}$ (4.42) to the interval $[0,\pi]$ which is
{\it submersive} over the open interval $(0,\pi)$; that is, all its fibres
are diffeomorphic.

 \item This map $\wave{\al_I}$ (4.42) factors as
 $$
 \wave{\al_I}\colon\wave{\orGr(3,TX_p)}@>\wave{\det_I}>>
 \proj^1=S^2@>h>>[0,\pi],
 \tag5.33
 $$
where $\wave{\det_I}$ is the regularisation of the map (4.10) after
the blowup of $\LaGr$.
 \endroster

The main observation is 

 \proclaim{Proposition 5.6} The map $\wave{\det_I}$ (5.33) 
 $$
\wave{\orGr(3, TX_p)}\to \proj^1=S^2
 \tag5.34
 $$
is submersive everywhere over the target sphere;

2) All fibres of this map are diffeomorphic and
 $$
\wave{\det_I}\1(z)=\al_I\1(0).
 \tag5.35
 $$
 \endproclaim

In particular, as in the 2-dimensional case, the space
$\wave{\det_I}\1(e^{i\fie})$ of {\it super special Lagrangian subspaces}
is isomorphic to the space of 3/2-pseudoholomorphic subspaces:
 $$
\wave{\det_I}\1(e^{i\fie})=\al_I\1(0).
 \tag5.36
 $$

Moreover, as in the 2-dimensional case, for every fibre
$\wave{\al}_I\1(z),z\ne0,\infty$ there exists an extra phase map
 $$
\al_I\1(z)\to S^1,
 \tag5.37
 $$
which for the oriented Lagrangian Grassmannian coincides with the standard
phase map
 $$
\det \colon\LaGr(3, TX_p)\to S^1_p
 \tag5.38
 $$
described in \S1.

We now globalise these maps. First of all, one has the oriented
Grassmannisation of the tangent bundle:
 $$
\pi \colon\orGr(3, TX)\to X
 \tag5.39
 $$
which fibrewise is $\orGr(3, TX_p)$, and the fibration
 $$
 \wave{\pi} \colon \wave{\orGr(3, TX)}\to X
 \tag5.40
 $$
of blown up Grassmannians
 $$
\si \colon \wave{\orGr(3, TX)}\to\orGr(3, TX)
 $$
along $\LaGr(TX)$.

Secondly, we have to globalise the target 2-spheres. The geometric meaning
of formula (4.13) is that the target space of the globalising of the
fibre-by-fibre maps $\det_{I,w}$ (4.10) is the projective $\proj^1$-bundle
 $$
\proj(\Oh_X \oplus \Oh_X(-K_X)).
 \tag5.41
 $$
This vector bundle $\Oh_X \oplus \Oh_X(-K_X)$ isn't simple. In particular
its automorphism group contains the multiplicative group
 $$
\C^*\subset\Aut (\Oh_X \oplus \Oh_X(-K_X))
 \tag5.42
 $$
which preserves the decomposition. Of course we can consider this subgroup as the 
subgroup of the automorphism group of the line bundle $L_{-K}=\Oh_X(-K_X)$:
 $$
\C^*\subset\Aut (L_{-K}).
 \tag5.43
 $$

\definition{Definition 5.1} The choice of an isomorphism
 $$
 \bigwedge^2(\Oh_X \oplus \Oh_X(-K_X))=L_{-K}
 \tag5.44
 $$
is called a {\it global complex orientation} of $X$.
\enddefinition

Thus the space of global complex orientations is again $\C^*$.

 \remark{Remark} Actually, we will consider slightly deformed Algebraic
Geometry. Thus the initial almost complex structure on $X$ is a holomorphic
structure on $X$, and we have the isomorphism
 $$
\C^*=Aut (L_{-K})
 \tag5.45
 $$
of the group of holomorphic isomorphisms of the holomorphic line bundle
$L_{-K}$.

In the general case, this identification is defined up to any gauge
transformation of $L_{-K}$.
 \endremark

 \remark{Example. The CY case} In this case a global complex orientation is
given by a choice of a holomorphic 3-form $\theta$.
 \endremark

We now fix a global complex orientation $w$ of $X$ as in (5.12); then we
have the globalisation of the map (4.10) 
 $$
\det_{I} \colon\orGr(3, TX)\to \proj(\Oh_X \oplus \Oh_X(-K_X))
 \tag5.46
 $$
of the fibrations over $X$.

The restriction of this map to the oriented Lagrangian Grassmannisation $\LaGr (TX)$ 
gives the map 
 $$
 \det\colon\LaGr (TX)\to S^1 (L_{-K})\subset \proj(\Oh_X
\oplus
 \Oh_X(-K_X))
 \tag5.47
 $$
(just as in (2.42)).

This target bundle (a ruled fourfold) has two sections
 $$
S_0=\proj(\Oh_X) \quad \text{and} \quad S_\infty=\proj(\Oh_X(-K_X)),
 \tag5.48
 $$
and an $S^1$-subbundle
 $$
S^1(L_{-K})\subset \proj(\Oh_X \oplus \Oh_X(-K_X)).
 \tag5.49
 $$
and these three subspaces are pairwise disjoint.

The regularisation of the map (5.46) is the map
 $$
\wave{\det_{I}}\colon\wave{\orGr(3, TX)}\to \proj(\Oh_X \oplus
\Oh_X(-K_X)),
 \tag5.50
 $$
which is submersive and smooth.

Finally, over every point $p\in X$ the space of pairs $(I,w)$ of {\it
oriented complex structures} is a lifting of the projective space
$\proj^3_+$ to a vector space $W^+$:
 $$
W^+_p=\{(I, w)\}.
 \tag5.51
 $$

 \proclaim{Proposition 5.7} The space of oriented complex structures is
the vector bundle 
 $$
W^+=\Oh_X \oplus TX.
 \tag5.52
 $$
 \endproclaim 

The automorphism group $\Aut W^+$ containes $\C^*$ acting by homotheties on
$TX$. This subgroup can again be identified to $\C^*\subset\Aut L_{-K}$. 

The fibrewise blowdown of the subbundle $\proj TX$ to point gives the Thom
bundle
 $$
S=TX \cup X_\infty,
 \tag5.53
 $$
which fibre-by-fibre is a $S^6$-bundle (see (5.22)).

\head \S6. SpLag and sdAG 3-cycles \endhead

Now for any oriented 3-dimensional submanifold $Y\subset X$ the embedding
$i\colon Y\to X$ has a Gauss lifting:
 $$
G(i)\colon Y\to\orGr(3, TX) \vert_Y
 \tag6.1
 $$
sending a point $y\in Y$ to the oriented subspace $TY_y\subset TX_y$.
The composite of this Gauss map and the map $\det_I$ (5.46) gives the map
 $$
\det_I\colon Y\to \proj(\Oh_X \oplus \Oh_X(-K_X)) \vert_Y.
 \tag6.2
 $$

\definition{Definition 6.1} 1) An oriented 3-dimensional cycle $Y\subset
X$ is 2/3-{\it pseudo\-holo\-morphic} (2/3-ph for short) if 
 $$
\det_I (Y)\subset S_0
 \tag6.3
 $$
(see (5.48));

2) it is called 2/3-{\it anti-pseudoholomorphic} (2/3-aph for short) if
 $$
\det_I (Y)\subset S_\infty.
 \tag6.4
 $$
\enddefinition

It is easy to see that $Y$ is Lagrangian iff 
 $$
\det_I (Y)\subset S^1(L_{-K}).
 \tag6.5
 $$
As before, let $a_{L-C}$ be the Levi--Civita connection on the
anticanonical line bundle $L_{-K}$ and
 $$
F_{L-C}\in\Om^2 (X)
 \tag6.6
 $$
the curvature of this connection as a 2-form on $X$.

 \definition{Definition 6.2} An oriented 3-cycle $Y$ on $X$ is {\it
canonically flat} if there exists a simple connected submanifold $B\supset
Y$ such that
 $$
F_{L-C} \vert_Y=0.
 \tag6.7
 $$
 \enddefinition

As before, if $Y$ is canonically flat then there exists a canonical
trivialisation of the restriction of the target bundle of the map $\det_I$
 $$
 \proj(\Oh_X \oplus \Oh_X(-K_X)) \vert_Y=Y \times \proj^1 
 $$
and a canonical projection
 $$
\pr\colon \proj(\Oh_X \oplus \Oh_X(-K_X)) \vert_Y\to \proj^1=S^2
 \tag6.8
 $$

The composite of these maps gives the map of a canonically flat cycle $Y$:
 $$
m_I=\pr \cdot \det_\C \cdot G(i)\colon Y\to S^2.
 \tag6.9
 $$

To use the regularisation of this map, that is, the map (5.50), we have to
correct slightly the definition of 3-dimensional cycle: every 3-cycle $Y$
has the (compact) subset
 $$
Y_\La=\{y\in Y \bigm| TY_y\in\LaGr(3, TX_y)\},
 \tag6.10
 $$
that is, the set of points with Lagrangian tangent space.

Now on $Y \setminus Y_\La$ there is defined the field of directions
 $$
0\to\ker (\om \vert_Y)\to TY @>\om >> T^*Y
 \tag6.11
 $$
 
 \definition{Definition 6.3} A pair $(v, Y)$ where $v$ is any smooth
extension of the field of directions $\ker\om \vert_Y$ from $Y \setminus
Y_\La$ to $Y$ is called a {\it supercycle} if this field is parallel with
respect to the Levi--Civita connection on $TY$.
 \enddefinition

For example, any 3/2-ph cycle can have only one supercycle structure:
$Y_\La=\emptyset$ and $\om \vert_Y$ has to be parallel with respect to $g
\vert_Y$. We call it s3/2-ph cycle for short and we have to give the

 \proclaim{Warning} Not every 3/2-ph cycle is s3/2-ph cycle.
 \endproclaim

But for a Lagrangian cycle supercycle structure is given by any character
of the fundamental group of $Y$.

Now for any oriented super submanifold $(v,Y)$ there exists a Gauss
lifting of the embedding $i\colon Y\to X$:
 $$
 G(i)\colon (v,Y)\to \wave{\orGr(3, TX)} \vert_Y
 \tag6.12
 $$
sending a point $y\in Y$ to the oriented pair 
 $$
 v_yTY_y\in \wave{\orGr(3, TX_y)}.
 \tag6.13
 $$
The composite of this Gauss map and the map $\wave{\det_I}$ (5.50) gives
the map
 $$
 \det_\C\colon Y\to \proj(\Oh_X \oplus \Oh_X(-K_X)) \vert_Y.
 \tag6.14
 $$

\definition{Definition 6.5} 1) An oriented 3-dimensional supercycle
$(v,Y)$ is called a {\it Lagrang\-ian supercycle} if 
 $$
 \det_\C (Y)\subset S^1(L_{-K})
 \tag6.15
 $$
(see (5.49));
\enddefinition

So for any canonical flat 3-dimensional supercycle $(v,Y)$ the {\it
complex phase map}
 $$
m_\C={\pr} \circ {\det_\C} \circ G(i)\colon Y\to S^2.
 \tag6.16
 $$
is well defined as the composite of standard regularised maps. 

The target sphere $S^2$ of this map has two special points
 $$
 0=\proj^1 \cap S_0 \quad\text{and}\quad \infty=\proj^1 \cap S_\infty,
 \tag6.17
 $$
and the circle
 $$
 S^1=\proj^1 \cap S^1(L_K).
 \tag6.18
 $$
This sphere admits the standard complex structure, and one has the
decomposition (2.37--2.38)
 $$
 S^2=\C \cup\infty=D^+ \cup S^1 \cup D^-.
 $$
We can extend the list (6.3--6.5) by the following cases:
 $$
 \align
 m_\C (Y)\subset D^+ &\iff Y \text{ is 2/3-symplectic},\\
 m_\C (Y)\subset D^- &\iff Y \text{ is 2/3-antisymplectic},
 \endalign
 $$
and so on.

\definition{Definition 6.4} 
1) A canonically flat cycle $Y$ is a {\it sdAG cycle} if $m_\C (Y)$ is
a point in the target sphere $S^2$ (6.16), that is, the differential of
this map vanishes:
 $$
 \dd m_\C=0.
 \tag6.19
 $$

3) a sdAG cycle $Y$ is a {\it spLag cycle} if 
 $$
 m_\C (Y)\in S^1\subset S^2
 \tag6.20
 $$
(see (6.18); it is easy to see that this map forgets super structure);

4) a sdAG cycle $Y$ is called $\al$-cycle if 
 $$
 \al_I (Y)=\al\in [0,\pi].
 \tag6.21
 $$
\enddefinition 

It is easy to see that 

1) a sdAG cycle $Y$ is 2/3-ph iff
 $$
 m_\C (Y)=0\in S^2,
 \tag6.22
 $$
and is 2/3-aph iff
 $$
m_\C (Y)=\infty
 \tag6.23
 $$
(see Definition 5.1).

\subheading{Mirror digression: CY threefolds} Again, every Lagrangian
3-cycle is canonically flat if $X$ is a simply connected Calabi--Yau
threefold. In this case $B=X$ again and {\it there exists a global
complex phase map}
 $$
m_\C \colon \wave{\orGr(3, TX)}\to S^2
 \tag6.24
 $$
such that for every super 3-cycle $Y$ its complex phase map is the
composite of the Gauss map and this universal map.

Moreover every sdAG cycle $Y$ on a CY threefold $X$ defines a complex
orientation of $X$, that is, a trivialisation of the canonical line bundle
$\bigwedge^3 TX$. Such a trivialisation is given by a choice of a
holomorphic 3-form $\theta$. A pair $(X,\theta)$ is called an {\it oriented
CY threefold}.

A spLag cycle $Y$ on an oriented CY threefold $X$ is Lagrangian with
respect to the K\"ahler form $\om$ and satisfying the condition
 $$
 \Re \theta \vert_Y=0
 \tag6.25
 $$
(see [H-L]).

Our aim is to investigate the local deformation theory of super sdAG
cycles in complex threefolds just as we did in \S3 for the case of aK
surfaces (see (3.24--3.28)). Recall that the starting point of this
investigation was the theory of complete linear systems of holomorphic
curves on algebraic surfaces. The main fact about such deformation theory
was Observation~3.1. In 3-dimensional case our experience is the local
deformation theory of spLag cycles on CY threefolds.

Such local deformation theory is quite well understood (for a survey, see
for example [H]). Let $\sM_Y$ be the local deformation space of spLag
cycles around a smooth spLag cycle $Y$. Then the tangent space to  the
moduli space at the point $Y$ is
 $$
T \sM_Y=H^1(Y,\R)
 \tag6.26
 $$
the space of harmonic 1-forms (in the induced metric) on $Y$. 

This space is the space of infinitesimal deformations and the obstructions
space for the Kuranishi family of deformations is $H^2(Y,\R)$. However,
onece more, there are no genuine obstructions: every infinitesimal
deformation extends to a geometric deformation (just as for holomorphic
curves on an algebraic surface $S$ of positive geometric genus, see
(3.17)).

On the other hand, if $Y$ is a homology sphere (that is, $H^*(Y)=H^*(S^3)$)
then this spLag cycle is {\it rigid}.

Now let $S\sM_{(v,Y)}$ be the local deformation space of super spLag
cycles around a smooth super spLag cycle $(v,Y)$. Then the tangent space
to the moduli space at the point $(v,Y)$ is
 $$
T S\sM_{(v,Y)}=H^1(Y,\R)\oplus i \cdot H^1(Y,\R)=H^1(Y,\C)
 \tag6.27
 $$
as the space of harmonic {\it complex} 1-forms (in the induced metric) on
$Y$. 

This space is the space of infinitesimal deformations and the obstructions
space for the Kuranishi family of deformations is $H^2(Y,\C) $. However
there are no genuine obstructions  again: every infinitesimal deformation
extends to a geometric deformation.

But now, just as in the 2-dimensional case
 $$
 \align
 &\text{the Geometry of special Lagrangian supercycles can be deformed}\\
 &\text{to the Geometry of sdAG 2/3-pseudoholomorphic supercycles}.
 \endalign
 $$
Indeed, let $e^{i\fie}$ be the image of spLag supercycles. Then the space
of special super Lagrangian directions is
 $$
m_\C\1(e^{i\fie})\subset S^2(U)\subset \wave{\orGr(3, TX)},
 \tag6.28
 $$
where $m_\C$ is the universal complex phase map (2.24), $S^2(U)$ is the
unit sphere bundle of the tautological bundle over the Lagrangian
Grassmannisation $\LaGr(3,TX)$ of the tangent bundle to $X$.

Now we can deform this fibre along the meridian to the north pole 0 of the
target space of the universal complex phase map (2.24) (see (6.17)). So we
get a family of geometries parametrised by the interval $[0,\pi/2]$ 
 $$
\Cal G=m_\C\1([0,\pi/2])
 \tag6.29
 $$
which gives a {\it smooth bordism} between the {\it special Lagrangian
super Geometry} given by the space of super directions (6.28) and the
{slightly deformed Algebraic Geometry} given by the space of
3/2-pseudoholomorphic directions
 $$
m_\C\1(0)\in\orGr(3, TX).
 \tag6.30
 $$
This is the explanation why we have the local deformation theory for
s3/2-ph cycles as good as the the local deformation theory for spLag
supercycles: let $S3/2\sM_Y$ be the local deformation space of
s3/2-ph cycles around a smooth s3/2-ph cycle $Y$. Then 

 \proclaim{Proposition 6.1} 1) The tangent space to 
the moduli space at the point $Y$ is
 $$
T S3/2\sM_Y=H^1(Y,\C);
 \tag6.31
 $$

 2) this space has a ``canonical'' complex orientation;

 3) there are no genuine obstructions: every infinitesimal deformation
extends to a geometric deformation.
 \endproclaim 

For the proof we have to imitate the arguments of Hitchin [H] and McLean in
the complex version.

\subheading{Remark} We should remark that our proposed complex version of
the theory of spLag cycles is not contained within symplectic geometry.
Indeed, sdAG cycles don't have to be symplectic in general. We consider the
metric as the fundamental object, not the symplectic form $\om$. Indeed if
we are changing the almost complex structure $I_t$ starting with an aK
triple $(\om_0,I_0, g)$ preserving our metric $g$ and such that $\om_t$
are harmonic, then, for small enough $t$, in the Hermitian triple
$(\om_t,I_t,g)$ the 2-form $\om_t$ is positive, but after some time it
loses this property and we go out of symplectic stuff.

One of the reasons why it is fruitful to go outside symplectic geometry is
the main Hitchin construction from [H]. The point is that locally the
moduli space $\sM_Y$ of deformations of spLag cycles in fixed CY threefold
$X$ around $Y$ can be embedded in the space $H^1(Y,\R) \times H^2(Y,\R)$
and the image is spLag cycle with respect to the standard symplectic
structure (recall that $H^1(Y,\R)=H^2(Y,\R)^*$), the standard complex
structure which define the metric. This metric isn't Riemannian but
ultrahyperbolic (but the restriction to $Y$ is Riemannian). See [H]. It is
easy to see that one can change this triple canonically in such a way that
the restriction to $Y$ doesn't change and the new metric is Riemannian but
the image of $\sM_Y$ isn't symplectic, and with respect to this new triple
the image of the moduli space is sdAG but not spLag.

So the complex version of the spLag cycles is productive to preserve the
standard pattern: ``a moduli space of complex submanifolds is a complex
manifold'' (like complete linear systems on an algebraic surface).

 \proclaim{Observation} The moduli space of sdAG supercycles is sdAG.
 \endproclaim

The last remark in the Calabi--Yau set-up is the following: if you believe
in Strominger, Yau and Zaslov's version of mirror symmetry (see [SYZ]), we
can expect the following:

 \proclaim{Great Expectations \rm [D]} 1) There exists a natural
compactification
 $$
 \vert Y \vert
 \tag6.32
 $$
of the global moduli space of all s3/2-ph deformations of a smooth
s3/2-ph cycle $Y$ in a Calabi--Yau threefold $X$ (we call it a {\rm
complete linear system} of s3/2-ph cycles).

2) The natural complex structure (given by the equality (6.31) is
integrable and extends to the complete variety $\vert Y \vert$; then
 $$
 Y=T^3=S^1 \times S^1 \times S^1 \implies
 \vert Y \vert=X' \quad \text{is a CY threefold}.
 \tag6.33
 $$
 \endproclaim

But even outside the Calabi--Yau set-up, these constructions give a
beautiful family of Geometries. We conclude our collection of definitions
and constructions by extending them to complex threefolds with positive and
negative canonical class.

In the set-up of CY threefolds, one can reproduce the following collection
of notions: suppose $Y$ is a 2/3-symplectic oriented cycle in $X$. Then the
image $m_\C (Y)\subset D^+$ is a compact subset and there exists a unique
minimal disc containing this image
 $$
m_\C (Y)\subset D_Y\subset D^+ 
 \tag6.34
 $$
(see (2.39)).

Again,
 \roster
 \item the centre $c_Y$ of the disc $D_Y$ is called the {\it centre} of
$Y$;

 \item the radius $r_Y$ of the disc $D_Y$ is called the {\it radius} of
$Y$;

 \item $Y$ is called $\ep$-spLag if
 $$
 c_Y\in S^1 \quad \text{and} \quad r_Y < \ep
 \tag6.35
 $$
 \endroster
It would be great to get some analogue of Donaldson's Theorem 2.1 for 
$\ep$-spLag cycles.

Now consider an simple connected smooth algebraic threefold $X$ with
canonical class $K_X>0$ (or $K_X<0$). We would like to deform slightly the
2/3-ph geometry, by considering a distinguished family of oriented
3-cycles. Again these cycles are determined by their first order
infinitesimal behavior, that is, the cycles defined by properties of their
Gauss lifts.

We can lift the projective bundle $\proj (\Oh_X \oplus \Oh_X(-K_X))$ to a
vector bundle 
 $$
 \aligned
V_K=&\Oh_X \oplus \Oh_X(K_X) \quad \text{ if } K_X >0, \\
 \text{and} \quad
 &\Oh_X \oplus \Oh_X(-K_X) \quad \text{ if } K_X < 0.
 \endaligned
 \tag6.36
 $$

\subheading{General type \rm($K_S>0$)} In this case, we consider a
nonvanishing section $s$ (or a section vanishing along a ``divisor'') 
to get a section of the projective bundle $\proj(\Oh_X\oplus\Oh_X(-K_X))$,
like $S_0$ and $S_\infty$ (5.48). 

In the K\"ahler (algebraic) case we have a finite dimensional family of
{\it holo\-morphic} sections
 $$
 s\in H^0(V_K)=\C \oplus H^0(\Oh_X(K_X))
 \tag6.37
 $$
We write $p_g=\dim H^0(\Oh_X(K_X))$ for the geometric genus of $X$, that
is, the complex dimension of $H^{3,0}(S)$. The family of such sections
defines a family of sections of the projectivisation $\proj V_K$ 
 $$
 S^{2p_g}=\C^{p_g} \cup \{\infty\setminus\text{point}\},
 \tag6.38
 $$
where $\C^{p_g}=\{(1, s)\}$ is the space of nonvanishing sections
and $\proj H^0(\Oh_X(K_X))=\vert K_X \vert$ is the complete canonical
linear system, points of which give the same linear subbundle, that is,
the same section of the projectivisation. That is, the map
 $$
\si_\infty\colon \proj H^0(\Oh_X(K_X))\to S^{2p_g}
 \tag6.39
 $$
isn't holomorphic, and blows the hyperplane $\vert K_X \vert$ down
to the point $\infty$.

Now the Levi--Civita connection gives the Hermitian structure on the
canonical bundle $L_K=\Oh_X(K_X)$ and, similarly, the Hermitian structure
on $\C^{2p_g} \cup \{\infty\setminus\text{point}\}$ gives the standard
metric on this sphere. Thus we can identify our sphere with the dual
sphere (see (2.57--2.60)).
 
So this sphere contains the ``equator''
 $$
 S^{2p_g-1}_e=\{z \bigm| \Vert z \Vert=1\}\subset\C^{2p_g}.
 \tag6.40
 $$
Again interpretation of $\C^{2p_g} \cup \{\infty\setminus\text{point}\}$
as the space of sections gives us the embedding
 $$
 i_{\can} \colon \proj V_K\to \proj H^0(V_K)^* \times X
 \tag6.41
 $$
and the projection of the trivial bundle to the fibre gives the map
 $$
 i_{\can} \colon \proj V_K\to \proj H^0(V_K)^*\to S^{2p_g}.
 \tag6.42
 $$
Finally, the composite of this projection and the blowdown (6.39) gives the
map
 $$
 \pr \colon \proj V_K\to S^{2p_g}.
 \tag6.43
 $$

Now for any cycle $Y\subset X$, the composite of the Gauss map,
the projection $\wave{\det_I}$ (5.50) and (6.43) defines the
complex phase map
 $$
 m_\C={\pr} \circ \wave{\det_I} \circ G(i) \colon X\to S^{2p_g}.
 \tag6.44
 $$
Now in terms of this phase map, one can define the analogues of the
sdAG and spLag cycles known in the Calabi--Yau case.

\definition{Definition 6.6} A cycle $Y\subset X$ is called a sdAG cycle
if $m_\C(Y)$ is a point, or equivalently
 $$
 \dd m_\C=0.
 \tag6.45
 $$
\enddefinition
\definition{Definition 6.7} 1) A 3-cycle $Y$ is called {\it weakly
Lagrangian} (wLag cycle for short) if
 $$
 m_\C(\Si)\subset S^{2p_g-1}_e,
 \tag6.46
 $$
where $S^{3p_g-1}_e$ is the equator ;

2) $Y$ is called a spLag cycle, if it is a sdAG cycle and
 $$
 m_\C(Y)\in S^{2p_g-1}_e.
 \tag6.47
 $$
\enddefinition

The equator divides the target sphere of the complex phase map into upper
and lower hemispheres:
 $$
S^{2p_g}=D^+ \cup S^{2p_g-1}_e \cup D^-,
 \tag6.48
 $$
and the entire catalogue of definitions (such as Definition 2.3),
properties and facts can be repeated in this new set-up.

\subheading{Fano case} Finally, it is quite easy to see what to do if
$K_S<0$: we change the sign of the canonical system $K_S\to-K_S$, getting
the sphere 
 $$
 S^{2h^0(\Oh_S(-K_S))}
 \tag6.49
 $$
as the target sphere of the complex phase map. After that, we can repeat
all our constructions and definitions.

\head \S7. Geometry of 3/2-pseudoholomorphic supercycles \endhead

Whereas the system of differential equations for spLag cycles are
notoriously complicated (see [H-L]), that for s3/2-ph cycles is much
simpler. But instead of describing it, we explain why s3/2
pseudoholomorphic Geometry is a slight deformation of Algebraic Geometry.
There is a type of complex oriented threefolds (that is, with trivial
canonical class) which are 2-connected, and thus can't be K\"ahler. In
this case the Geometry of s3/2-ph cycles is the unique tool for
investigations. It is proper to show how the s3/2 AG works in this
case.

Recall briefly the general properties of such type manifolds. The class of
such manifolds was introduced by Miles Reid [R] in the set-up of the
investigation of minimal models theory of complex threefolds. Almost at
the same time R. Friedman in [F] proposed the basic construction of such
manifolds in the set-up of his theory of infinitesimal deformations of
singular complex manifolds. So it is resonable to call these complex
manifolds {\it Friedman--Reid manifolds} (FR threefolds for short) in the
same vein as in the K\"ahler case, we call it Calabi--Yau (CY for short)
threefolds.

 \definition{Definition 7.1} A {\it FR threefold} is a compact smooth
complex threefold $X$ such that
 $$
H^p(X, \Om^q)=0 \quad \text{ for } p+q \ne0, 3, 6,
 \tag7.1
 $$
where $\Om$ is the cotangent bundle of $X$, and $X$ is $2$-connected:
 $$
\pi_n (X)=0 \quad \text{ for } n < 3.
 \tag7.2
 $$
\enddefinition

C.T.C. Wall proved that all FR threefolds are diffeomorphic to a connected
sum of $g$ copies of $S^3 \times S^3$, where $g$ is any positive integer,
which we call the {\it genus} of $X$. 

The main properties of FR threefolds are just the same as of CY threefolds:

 \proclaim{Main properties (\rm [C], Proposition 1.3)} For a FR threefold
$X$,
 \roster
 \item The Hodge spectral sequence 
 $$
 E_2^{p,q}=H^q(X, \Om^p)
 \tag7.3
 $$
degenerates at $E_2$. Thus the cohomology of $X$ has an integral Hodge
structure which is pure of weight 3.

 \item The local deformation space of $X$ is unobstructed, smooth with the
tangent space at $X$ equal to
 $$
 H^1(TX)=H^1(X,\Om^2)=\C^{g-1}.
 \tag7.4
 $$

 \item Let $F^2H^3(X,\C)$ be cohomology classes which can be represented by
d-closed forms of types (3,0) and (2,1). Then we get the orthogonal
decomposition
 $$
 H^3(X,\C)=F^2H^3(X,\C)\oplus \overline{F^2H^3(X,\C)}
 \tag7.5
 $$
 and the Hermitian form
 $$
(1/2i)\Span{\al, \overline{\be}}.
 \tag7.6
 $$
 \endroster
 \endproclaim

Recall that the signature of this Hermitian form (7.6) is called the {\it
signature} of $X$.

As usual we fix a {complex orientation} of $X$ that is a holomorphic 3-form
$\theta$ and the pair $(X,\theta)$ is called an oriented FR threefold. 

\subheading{The theory of periods of FR threefolds} In terms of the Hodge
structure on $X$, the existence of an integral 3-cohomology class
 $$
[Y]\in H^3(X,\Z)
 \tag7.7
 $$
which can be realised as a smooth 3/2-pseudoholomorphic cycle $Y$ gives one
equation on the period domain of FR threefolds. Indeed, as for K3
surfaces, in this case
 $$
 \int_Y \theta=0.
 \tag7.8
 $$

\definition{Definition 7.2} A 3-cohomology class $[Y]$ is called a {\it
polarisation} of $X$ if there exists a global moduli space $\vert Y
\vert$ of all deformations of $Y$ in $X$ as s3/2-ph cycles with the
(natural) complex structure (we call it a complete linear system) such
that 

1) the complete linear system $\vert Y \vert$ is an {\it algebraic variety} and 

2) the complex dimension
 $$
\dim \vert Y \vert > 2
 \tag7.9
 $$
\enddefinition

 \remark{Remark} You can see that the {\it dimension} of the moduli space
works in the same way as {\it Riemann--Roch formula} plus {\it adjunction
formula} in the case of complex surfaces.
 \endremark

Now we can define the analogue of the Picard lattice for FR threefolds:

\definition{Definition 7.3} The sublattice
 $$
3/2\text{-}\Pic(X)\subset H^3(X,\Z)
 \tag7.10
 $$
generated by s3/2-ph cycles is called the {\it $3/2$-Picard lattice} of
$X$.
\enddefinition

In the same vein, using the Hermitian form (7.6), we can define the
sublattice of {\it transcendental} cycles:
 $$
 \Tr(X)\subset H^3(X,\Z),
 \tag7.11
 $$
and these constructions are absolutely parallel to the standard
constructions for K3 surfaces.

We now return to the local geometry around a smooth s3/2-ph cycle $Y$.

For any 2/3-ph smooth cycle $Y$ we have the orthogonal ``Hodge
decomposition''
 $$
TY=\ker(\om \vert_Y)\oplus TY_I
 \tag7.12
 $$
where $TY_I$ is the complex directions of the tangent bundle of $Y$.

 \definition{Definition 7.4} A smooth s3/2-ph cycle $Y$ is called {\it
wrapped} iff there exists a smooth holomorphic curve
 $$
 C\subset Y
 \tag7.13
 $$
inside $Y$.
 \enddefinition
 
Now to define a wrapped cycle in terms of the holomorphic curve $C$, we
have to consider the normal bundle 
 $$
 N_{C\subset X} 
 \tag7.14
 $$
of our curve in $X$. The action of the complex structure operator $I$ on
the subspace
 $\ker(\om \vert_Y)$ (7.12) defines the subbundle
 $$
 L_Y=\Span{\ker(\om \vert_Y), I(\ker(\om \vert_Y))}.
 \tag7.15
 $$
which is a {\it complex} subbundle.
 
 \proclaim{Proposition 7.2} $L_Y$ is a holomorphic subbundle of the
holo\-morphic rank $2$ bundle $N_{C\subset X}$ (7.14).
 \endproclaim
 
 Roughly speaking our s3/2-ph cycle $Y$ wraps the zero section in
the line bundle $L_Y$.
 
 Now by the definition
 $$
 \text{3/2-ph cycle $Y$ is a supercycle} \implies \deg L_Y=0
 \tag7.16
 $$
 and we get
 
 \proclaim{Proposition 7.3} For s3/2-ph cycle $Y$,
 \roster
 \item the line bundle $L_Y$ admits a flat $\U(1)$-connection,
 
 \item that is, there exists the character of the fundamental group of $Y$
 $$
 \rho \colon\pi (Y)\to\U(1)
 \tag7.17
 $$
which gives our line bundle:
 $$
 L_{\rho}=L_Y,
 \tag7.18
 $$
and

 \item
 $$
 Y=S^1 (L_{\rho})
 \tag7.19
 $$
is the unit circle bundle of this flat $\U(1)$-bundle. In particular,
topologically
 $$
 Y=C \times S^1.
 \tag7.20
 $$
 \endroster
 \endproclaim
 
Let us return to FR threefolds. Such manifold has no nontrivial line
bundles, and has no divisors. But some FR threefolds have an infinity of
mutually disjoint rigid elliptic curves. Consider such FR threefold $X$
with such a set $\{C_i\}$ of holomorphic elliptic curves. The normal bundle
of any such curve is
 $$
 N_{C\subset X}=\Oh_C(\xi)\oplus\Oh_C(-\xi) \quad\text{for }
\xi\in\Pic_0(C),
 \tag7.21
 $$
with
 $$
 h^0(\Oh_C(\xi))=h^1(\Oh_C(\xi))=h^1(\Oh_C(-\xi))=h^0(\Oh_C(-\xi))=0
 \tag7.22
 $$
and any line subbundle $L$ of degree 0 is either $\Oh_C(\xi)$ or
$\Oh_C(-\xi)$. So in this case we have two wrapped s3/2-ph cycles $Y_+$ and
$Y_-$ which topologically are
 $$
 C \times S^1=T^3=S^1 \times S^1 \times S^1. 
 \tag7.23
 $$

Then modulo the Great Expectations we have the complex 3-dimensional
manifold  $\vert Y_+ \vert$ and the question is:

 \proclaim{Typical question} How many other elliptic wrapped s3/2-ph
cycles are in $\vert Y_+ \vert$?
 \endproclaim

 \remark{Remark} You can see yourself what happens if $N_{C\subset X}$ is
a twisted Atiyah bundle (we get one s3/2-ph cycle) or the trivial bundle
twisted by a point of second order of $Pic_0(C)$ (we get a $S^2$-family of
s3/2-ph cycles).
 \endremark

We finish with the following construction: 3/2-ph cycles in $X$ are what
we need for the theory of FR structures. Let $X$ be a FR threefold, 
 $$
 [Y]\in 3/2\text{-}\Pic (X)
 \tag7.24
 $$
(see (7.10)) and $Y_+$ and $Y_-$ are two s3/2ps cycles from this cohomology class.

\definition{Definition 7. 5} Such two cycles are called {\it K3 cobordant}
if there exists a smooth 4-submanifold $S$ which is a {\it complex
symplectic surface} in the sense of \S2 of [D2] such that 
 $$
 \partial S=Y_+ \cup Y_-;
 $$
that is, $S$ is a complex surface with a {\it holomorphic symplectic form} 
 $$
\theta_2\in \Om^{2,0}.
 \tag7.25
 $$
\enddefinition

Of course this construction relates closely to the ``complexification'' of
the diffeo\-morphism group of $Y_{\pm}$ and to Nahm's equations (see
[D2]) but on the other hand, this complex symplectic bordism gives a
special loop in the intermediate jacobian 
 $$
 J^3(X)=F^2H^3(X,\C)/H^3(X,\Z).
 \tag7.26
 $$
This set-up is closely related to the theory of vector bundles on FR and
CY threefolds, and you can consider the paper [T] as the continuation of
these considerations.

On the other hand if $Y_{\pm}$ are wrapped on curves $C_{\pm}$ the rank 2
vector bundle $E$ on the complex symplectic bordism $S$ such that
 $$
 E \vert_{C_{\pm}}=N_{C_{\pm}\subset X}
 \tag7.27
 $$
we can consider as the analogue of a rational function between two
divisors. What a large new field of Numerical Geometry! So, you can see
that this subject is open for investigations and You can move this
Geometry Yourself. Good Luck!

\head Acknowlegments \endhead

This is the version of my talks at the March 1998 COW weekend meeting and
the Symplectic Geometry Workshop in Warwick University. I would like to
express my gratitude to Mathematical Institute of Warwick University and
personally to Miles Reid, Dietmar Salamon, Victor Pidstrigach and Izette
Weiss-Pidstrigach for support and hospitality. I wish to thank S.~Mukai
for many helpful suggestions and especially and H. Clemens for permission
to use his unpublished paper [C].

\enddocument